\documentclass[11pt,reqno]{amsart}
\usepackage[colorlinks=true,
pdfstartview=FitV, linkcolor=cyan, citecolor=magenta,
urlcolor=blue]{hyperref}
\usepackage{amsmath,amsfonts,latexsym,amssymb}
\usepackage{mathrsfs}
\usepackage[latin1]{inputenc}
\usepackage[T1]{fontenc}
\usepackage{ae,aecompl}
\usepackage{braket}
\usepackage{comment}
\usepackage{color}
\usepackage{graphicx}
\usepackage{subfig}
\usepackage{bbm}

\usepackage{enumitem}

\newtheorem{theorem}{Theorem}[section]

\newtheorem{corollary}[theorem]{Corollary}

\newtheorem*{thm*}{\protect\theoremname}
\theoremstyle{definition}

\newtheorem*{cor*}{\protect\corollaryname}

\long\def\@savemarbox#1#2{\global\setbox#1\vtop{\hsize\marginparwidth 
  \@parboxrestore\tiny\raggedright #2}}
\usepackage{mathrsfs }

\providecommand{\corollaryname}{Corollary}
\providecommand{\theoremname}{Theorem}

\begin{document}

\title{Kleinian viewpoints on higher rank worlds}
\author[Canary]{Richard Canary}
\address{University of Michigan}
\thanks{Canary was partially supported by grant DMS-2304636  from the National Science Foundation.}

\begin{abstract} 
This talk is designed to attract people who work on real hyperbolic manifolds to consider thinking about discrete subgroups
of higher rank Lie groups. To that end, we breezily discuss some applications of the ideas from the theory of Kleinian groups in the higher rank setting.
\end{abstract} 

\maketitle

\tableofcontents

\section{Introduction}

In this talk, I will highlight some recent applications of ideas and techniques from the theory of Kleinian groups
to the study of discrete subgroups of higher rank Lie groups.
I will begin with a brief introduction for those who, like me, are not totally conversant in the general theory of Lie groups.
I hope this will attract others working in the rank one setting to consider studying discrete groups in higher rank. A longer attempt to make the theory
of discrete subgroup of higher rank Lie groups accessible to rank one people is  my {\em Informal Lecture Notes on Anosov Representations}
which are available on my webpage and at the AMS Notes pages.

Since my choice of topics reflect my idiosyncratic personal interests, I will recommend some other surveys. The surveys of
Kapovich, Leeb and Porti \cite{KLP-survey,KL-survey} explain their approach to discrete subgroups which focuses on the action
on the quotient symmetric space and explores many parallels with the rank one setting. The surveys of Wienhard \cite{wienhard-icm}  and 
Burger-Iozzi-Wienhard \cite{BIW-survey} do an excellent job of presenting results motivated by the techniques  and ideas from  the Fuchsian setting.
Kassel's survey articles \cite{kassel-icm,kassel-survey} contain nice discussions of the relationship with projective geometry, 
see also Marquis \cite{marquis-survey}.
Finally, Canary-Zhang-Zimmer \cite{CZZ-survey} survey recent developments in Patterson-Sullivan theory in higher rank
(see alse Sambarino \cite{sambarino-dichotomy}, which is not a survey paper but contains a very clear discussion of the Anosov case).

\medskip\noindent
{\bf Acknowledgements:} I would like to thank Martin Bridgeman,  Subhadip Dey, James Farre, Sourav Ghosh, Andres Sambarino, Kostas Tsouvalas, Gabriele Viaggi, Teddy Weisman, 
Neza $\check{\mathrm{Z}}$ager Korenjak, Tengren Zhang and Andy Zimmer
for helpful comments on an early version of this paper.

\section{Types of discrete subgroups}

\subsection{Linear algebra}
We will work mostly in the setting of discrete subgroup of $\mathsf{PSL}(d,\mathbb K)$ where 
$\mathbb K$ is $\mathbb R$ or $\mathbb C$. In this setting, much of the Lie theory is simple linear
algebra involving singular values and eigenvalues.
The Cartan subspace $\mathfrak a$ of
$\mathsf{PSL}(d,\mathbb K)$ is the space of real
trace-free diagonal matrices, which me may think of as $d$-tuples of real numbers whose entries add up to 0.
The (closed) positive Weyl chamber $\mathfrak a^+$ is simply $d$-tuples in $\mathfrak a$ whose entries are in descending order, i,e.
$$\mathfrak a=\{\vec a\in \mathbb R^d:a_1+\cdots+a_d=0\}\quad\text{and}\quad
\mathfrak a^+=\{\vec a\in\mathfrak a:a_1\ge\cdots\ge a_d\}.$$

We define the  Jordan projection $J:\mathsf{PSL}(d,\mathbb K)\to \mathfrak a^+$ and the
Cartan projection $\kappa:\mathsf{PSL}(d,\mathbb K)\to \mathfrak a^+$  to be
$$J(\gamma)=(\log \lambda_1(\gamma),\ldots,\log\lambda_d(\gamma))\quad\text{and}\quad
\kappa(\gamma)=(\log \sigma_1(\gamma),\ldots,\log\sigma_d(\gamma))$$
where $\lambda_i(\gamma)$ is the modulus of the $i^{th}$ eigenvalue of $\gamma$ and
$\sigma_i(\gamma)$ is the $i^{th}$ singular value of $\gamma$.

Let $X_d=\mathsf{PSL}(d,\mathbb R)/\mathsf{PSO}(d)$ (or $X_d=\mathsf{PSL}(d,\mathbb C)/\mathsf{PU}(d)$) be the symmetric space of 
$\mathsf{PSL}(d,\mathbb R)$ (or $\mathsf{PSL}(d,\mathbb C)$).
One should think of $J(A)$ and $\kappa(A)$ as vector-valued distance functions on $X_d$. To make this precise, note that
if  $x_0=[\mathsf{PSO}(d)]\in X_d$ (or $x_0=[\mathsf{PU}(d)]\in X_d$ if $\mathbb K=\mathbb C$), then 
$$d_{X_d}(x_0,A(x_0))=||\kappa(A)||\quad\text{and}\quad \inf_{x\in X_d}d_{X_d}(x,A(x)))=||J(A)||.$$
(Your favorite metric may be some constant multiple of our metric, but this choice will be convenient for us. For example, with this
convention, $X_2=\mathbb H^2$ with constant curvature $-\frac{1}{4}$.)

We may write any element $A\in\mathsf{PSL}(d,\mathbb K)$ as $Ke^{\kappa(\gamma)}L$ where $K,L\in \mathsf{PO}(d)$
(if $\mathbb K=\mathbb R$)
or $K,L\in \mathsf{PU}(d)$ (if $\mathbb K=\mathbb C$). (This is known as the Cartan decomposition.)
The elements $K$ and $L$ are not uniquely defined, but 
if $\sigma_k(\gamma)>\sigma_{k+1}(\gamma)$,
then
$$U_k(\gamma)=K(\langle e_1,\ldots, e_k\rangle)\quad\text{where }\{e_1,\ldots,e_d\}\text{ is the standard basis for }\mathbb K^d$$ 
is well-defined and is the image of the $k$-plane which is stretched the most by $\gamma$.

We will also be interested in linear functionals in $\mathfrak a^*$. Particularly important are the simple roots $\alpha_k:\mathfrak a\to\mathbb R$
and the fundamental weights  $\omega_k:\mathfrak a\to\mathbb R$ given by
$$\alpha_k(\vec a)=a_k-a_{k+1}\quad\text{and}\quad \omega_k(\vec a)=a_1+\cdots+a_k.$$

A (full) flag in $\mathbb K^d$ is a collection $\{F^k\}_{k=1}^{d-1}$ where each $F^k$ is a $k$-dimensional $\mathbb K$-subspace of $\mathbb K^d$
and $F^j\subset F^k$ if $j<k$.
If $\mathcal F$ is the space of full flags in $\mathbb K^d$, we can define the Iwasawa cocycle 
$B:\mathsf{PSL}(d,\mathbb R)\times\mathcal F\to \mathfrak a$
by the relation
$$\omega_k(B(A,F))=\frac{||A(v_1)\wedge\cdots\wedge A(v_k)||}{||v_1\wedge\cdots\wedge v_k||}\quad\text{where}\quad
\{v_1,\ldots,v_k\} \text{ is a basis for }F^k\quad\text{for all }k$$
where $F^k$ is the $k$-dimensional subspace of $F$.
In rank one, the Iwasawa cocycle is equivalent to the Busemann cocycle.

If $\theta\subset\{1,\ldots,d-1\}$, a $\theta$-flag is a collection $\{F^k\}_{k\in\theta}$ where each $F^k$ is a $k$-dimensional 
$\mathbb K$-subspace of $\mathbb K^d$
and $F^j\subset F^k$ if $j <k$ and $j,k\in\theta$. Let
$\mathcal F_\theta$  be the space of $\theta$-flags. We then define a partial Cartan subspace
$$\mathfrak a_\theta=\{\vec a\in\mathfrak a :\alpha_j(\vec a)=0 \text{ if }j\notin \theta\}\quad\text{and}\quad
\mathfrak a_\theta^+=\mathfrak a_\theta\cap\mathfrak a^+.$$
Then the group $\mathfrak a_\theta^*$ of linear functionals  on $\mathfrak a_\theta$ has basis $\{\omega_k\}_{k\in\theta}$.
There is a projection $p_\theta:\mathfrak a\to\mathfrak a_\theta$ defined by the relation
$$\omega_k(\vec a)=\omega_k(p_\theta(\vec a))\quad\text{ for all }k\in\theta.$$
We can also define a partial Iwasawa cocycle just as before
$$\omega_k(B_\theta(A,F))=\frac{||A(v_1)\wedge\cdots\wedge A(v_k)||}{||v_1\wedge\cdots\wedge v_k||}\quad\text{where}\quad
\{v_1,\ldots,v_k\} \text{ is a basis for }F^k\quad\text{for all } k\in\theta.$$

\subsection{Divergent and transverse groups}
A group $\Gamma\subset \mathsf{PSL}(d,\mathbb K)$ is {\em $P_k$-divergent} if whenever $\{\gamma_n\}$ is a sequence
of distinct elements in $\Gamma$ we have $\alpha_k(\kappa(\gamma))\to \infty$ where 
$$\alpha_k(\vec a)=a_k-a_{k+1}$$
is the $k^{\rm th}$ simple root. Intuitively, a $P_k$-divergent group is a discrete group whose discreteness is detected by the $k^{\rm th}$
simple root. Notice that if $\Gamma$ is $P_k$-divergent, then it is $P_{d-k}$-divergent, since $\alpha_k(\gamma)=\alpha_{d-k}(\gamma^{-1})$.
If $\Gamma$ is $P_k$-divergent, then both $U_k(\gamma)$ and $U_{d-k}(\gamma)$ are well-defined for all but finitely many
$\gamma\in \Gamma$. 

We say that $\Gamma$ is $P_\theta$-divergent for some non-empty $\theta\subset\{1,\ldots,d-1\}$ if it is $P_k$-divergent for all $k\in\theta$.
We will assume from now on that $\theta$ is symmetric.
Then $U_\theta(\gamma)=\{U_k(\gamma)\}_{k\in\theta}$ is a well-defined element of $\mathcal F_\theta$ for all but 
finitely many $\gamma\in\Gamma$.
The limit set $\Lambda_{\theta}(\Gamma)$ of $\Gamma$ is then just the set of accumulation points of
$\{U_\theta(\gamma)\}_{\gamma\in\Gamma}$ in $\mathcal F_\theta$.

If $\phi\in\mathfrak a_\theta^*$, we may define a ``distance function'' on the orbit $\Gamma(x_0)$ in $X_d$ by $d_\phi(x_0,\gamma(x_0))=\phi(\kappa(\gamma))$,
although of course this ``distance'' may be negative. (Notice that in the constant curvature $-1$ metric on $\mathbb H^2$, we have
$d(x_0,\gamma(x_0))=\alpha_1(\kappa(\gamma))$.)
We may then define
$$\delta_\phi(\Gamma)=\limsup_{T\to\infty} \frac{\log\#\{\gamma\in\Gamma: \phi(\kappa(\gamma))\le T\} }{T},$$
which we regard as the exponential growth rate of the orbit from the point of view of the linear functional $\phi$.
This is also the critical exponent of the $\phi$-Poincar\'e series of $\Gamma$. If $\delta^\phi<+\infty$, then
one may define a $\phi$-Patterson-Sullivan
measure supported on $\Lambda_\theta(\Gamma)$ (see \cite{CZZ3}), i.e. a probability measure $\mu$ so that
$$\frac{d\gamma_*\mu}{d\mu}(F)=e^{-\delta_\phi(\Gamma) B_\theta(\gamma^{-1},F)}$$
for all $\gamma\in\Gamma$ and $F\in\Lambda_\theta(\Gamma)$.

A $P_\theta$-divergent group $\Gamma\subset \mathsf{PSL}(d,\mathbb K)$ is {\em $P_\theta$-transverse} if every two distinct flags
$F$ and $G$ in $\Lambda_\theta(\Gamma)$ are transverse, i.e. if $k\in \theta$, then $F^k\oplus G^{d-k}=\mathbb K^d$. One key feature of
transverse groups is that the action of $\Gamma$ on $\Lambda(\Gamma)$ is a convergence group action. Moreover,
one may establish a Hopf-Tsuji-Sullivan dichotomy for Patterson-Sullivan measures on $\Lambda_\theta(\Gamma)$ (see \cite{CZZ3} and
\cite{KOW}).

\subsection{Anosov groups}
A finitely generated group $\Gamma\subset \mathsf{PSL}(d,\mathbb K)$
 is said to be {\em $P_\theta$-Anosov} if there exists $a,C>0$ so that
$$\alpha_k(\kappa(\gamma))\ge a|\gamma|-C\text{ for all }\gamma\in\Gamma\ \  \text{ and }k\in\theta$$
where $|\gamma|$ is the word length of $\gamma$ with respect to some fixed generating set for $\Gamma$.
One may think of this as saying that $\Gamma$ is quasi-isometrically embedded  (with respect to the usual
bi-invariant metric on $\mathsf{PSL}(d,\mathbb K)$) and the fact
that it is quasi-isometrically embedded is detected by the $k^{\rm th}$-simple root for all $k\in\theta$.

 If $\Gamma$ is $P_\theta$-Anosov,
then it is Gromov/word hyperbolic and $P_\theta$-transverse and there exists a $\Gamma$-equivariant homeomorphism
$$\xi:\partial \Gamma\to \Lambda_{\theta}(\Gamma)$$
where $\partial \Gamma$ is the Gromov boundary of $\Gamma$.

Anosov groups have two important features. 1) The orbit map of $\Gamma$ into $X_d$ is a
quasi-isometric embedding, i.e. $d(x_0,\gamma(x_0))$ grows linearly in the word length of $\gamma$. 
2) There exists a neighborhood $U$ of the inclusion map in 
$\mathrm{Hom}(\Gamma,\mathsf{PSL}(d,\mathbb K))$ so that if $\rho\in U$, then $\rho(\Gamma)$ is $P_\theta$-Anosov.

In the case of $\mathsf{PSL}(2,\mathbb K)$, a subgroup is $P_1$-Anosov if and only if it is convex compact.
(Recall that in a rank one Lie group, a discrete group is convex cocompact if and only if its orbit map is a quasi-isometric
embedding.)

If $\Gamma\subset\mathsf{PSL}(d,\mathbb K)$ is $P_\theta$-Anosov, $\phi\in\mathfrak{a}_\theta^*$ and $\delta^\phi(\Gamma)<+\infty$,
then Sambarino \cite{sambarino-quant} (see also \cite{BCLS}) constructed a metric Anosov flow which is H\"older orbit equivalent to the geodesic
flow of $\Gamma$ and whose periods are exactly $\{\phi(J(\gamma))\}_{[\gamma]\in[\Gamma]}$. This flow is one of the main tools
in the study of  the dynamics of the action of $\Gamma$. It allows us to think of $\phi(J(\gamma))$ as giving a length function on the group.
(His construction can be generalized to the setting of transverse groups, see \cite{CZZ3}.)

\medskip\noindent
{\bf Historical remarks:}  The book of Benoist and Quint \cite{benoist-quint} is an excellent source for the Lie theory
discussed here. Divergent and transverse groups were first studied Kapovich, Leeb and Porti (see \cite{KLP}),
who called them regular and regular antipodal groups.  Albuquerque \cite{albuquerque} and Quint \cite{quint} were the
first to study Patterson-Sullivan measure for discrete subgroups of higher rank Lie groups. 

Labourie \cite{labourie-invent} first defined Anosov groups in his study of Hitchin representations. The general theory was further developed
by Guichard-Wienhard \cite{GW}, Gu\'eritaud-Guichard-Kassel-Wienhard \cite{GGKW}, Kapovich-Leeb-Porti \cite{KLP},
Bochi-Potrie-Sambarino \cite{BPS} and others. Labourie's original definition was dynamical. The equivalent definition we gave
is due to Kapovich-Leeb-Porti \cite{KLP-survey} and Bochi-Potrie-Sambarino \cite{BPS}.
Dey and Kapovich \cite{DK-PS} developed a theory of Patterson-Sullivan measures for Anosov groups.

Given a semi-simple Lie group $\mathsf{G}$ of non-compact type and a parabolic subgroup $P$ of $\mathsf{G}$,
Labourie \cite{labourie-invent} introduced a theory of $P$-Anosov subgroups of $\mathsf{G}$ which had the same
important features as above. Guichard and Wienhard \cite[Prop. 4.3, Remark 4.12]{GW} observed that  given $\mathsf{G}$ and $P$,
there exists an irreducible representation $\tau:\mathsf{G}\to\mathsf{PSL}(d,\mathbb R)$ (for some $d$) so
that $\Gamma\subset\mathsf{G}$ is $P$-Anosov if and only if $\tau(\Gamma)$ is $P_{1,d-1}$-Anosov. So many, but not all,
problems concerning Anosov groups can be studied by studying the $\mathsf{PSL}(d,\mathbb R)$ case and using this equivalence.

\section{Which groups can be Anosov?}

It is natural to ask which isomorphism classes of groups arise at Anosov groups. The most obvious examples
are the convex cocompact subgroups of rank one Lie groups. For example, convex cocompact subgroups of
$PO(d,1)=\mathrm{Isom}(\mathbb H^d)$ are $P_1$-Anosov subgroups of $\mathsf{PGL}(d+1,\mathbb R)$.

The next collection of examples, historically, is provided by  fundamental groups of
strictly convex, closed (real) projective manifolds. We recall that an open domain $\Omega\subset\mathbb P(\mathbb R^{d+1})$ is strictly convex
if it is a bounded strictly convex subset of some affine chart for $\mathbb P(\mathbb R^{d+1})$. If $\Gamma\subset \mathsf{PSL}(d+1,\mathbb R)$
preserves and acts properly discontinuously and  cocompactly on a strictly convex domain $\Omega\subset\mathbb P(\mathbb R^{d+1})$,
then $\Omega/\Gamma$ is a closed strictly convex projective manifold. Benoist \cite{benoist-divisible1} proved that in this case $\Gamma$ is $P_1$-Anosov.
However, $\Gamma$ is not $P_k$ Anosov for any $2\le k\le d-2$ (see \cite[Cor. 1.4]{canary-tsouvalas}).
Kapovich \cite{kapovich-gt} showed that certain Gromov-Thurston $d$-manifolds admit strictly convex projective structures so 
that their fundamental groups
arise as $P_1$-Anosov subgroups of $\mathsf{PSL}(d+1,\mathbb R)$.

Danciger, Gu\'eritaud, Kassel, Lee and Marquis \cite{DGKLM} showed that any hyperbolic Coxeter group is isomorphic to an Anosov group. 
The following is the best current result.

\begin{theorem}{\rm (Douba-Fl\'echelles-Weisman-Zhu \cite{DFWZ})}
If a hyperbolic group acts properly and cocompactly on a $CAT(0)$
cube complex, then it  is isomorphic to an Anosov subgroup of $\mathsf{PSL}(d,\mathbb R)$ for some $d$.
\end{theorem}

It is still unknown whether or not there is a linear hyperbolic group which is not isomorphic to an Anosov group. In fact, only recently has it been shown
that there exist linear hyperbolic groups which are not isomorphic to a convex cocompact subgroup of some rank one Lie group. Tholozan
and Tsouvalas \cite{tholozan-tsouvalas} gave examples of linear hyperbolic groups which are not isomorphic to a discrete subgroup of any 
rank one Lie group. Their examples are doubles of lattices in $\mathsf{Sp}(n,1)$ along cyclic subgroups.

On the other hand, one can place restriction on groups which admit certain flavors of Anosov-ness for a given group.
For example, any torsion-free $P_1$-Anosov subgroup of $\mathsf{SL}(4,\mathbb R)$ is isomorphic to a convex cocompact
subgroup of $\mathsf{PSL}(2,\mathbb C)$, see Canary-Tsouvalas \cite[Thm. 1.2]{canary-tsouvalas}. 
Moreover, one can obtain bounds on the cohomological dimension of Anosov groups.

\begin{theorem}{\rm (Canary-Tsouvalas  \cite{canary-tsouvalas})}
If $\Gamma\subset \mathsf{PSL}(d,\mathbb R)$ is  torsion-free and $P_k$-Anosov and $k\le \frac{d}{2}$,
then
\begin{enumerate}
\item
If  $(d,k)$ is  not $(2,1)$, $(4,2)$,
$(8,4)$ or $(16,8)$, then $\Gamma$ has cohomological dimension at most $d-k$.
\item
If  $(d,k)$ is  either $(2,1)$, $(4,2)$,
$(8,4)$ or $(16,8)$, then $\Gamma$ has cohomological dimension at most $d-k+1$.
\end{enumerate}
\end{theorem}

\noindent
{\em Sketch of proof:} 
If $d=2$, then $\Gamma$ is a Fuchsian group and the result is obvious and the case $d=3$ can be handled separately
(see \cite[Thm. 1.1]{canary-tsouvalas}).

If $d\ge 4$, we fix $x_0\in \partial\Gamma$ and a $(d-k+1)$-dimensional subspace $V$ containing
$\xi(x_0)^{(d-k)}$ and define an injective  map $F:\partial\Gamma-\{x_0\}\to \mathbb P(V-\xi(x_0)^{(k)})$ by letting $F(y)$ be
the line $\xi(y)^{(k)}\cap V$. Therefore, $\partial\Gamma$ has topological dimension at most $d-k$. However, Bestvina and
Mess \cite{bestvina-mess} showed that the cohomological dimension of a torsion-free word hyperbolic group is exactly one
more than the topological dimension of its boundary. Therefore, in all cases $\Gamma$ has cohomological dimension at most $d-k+1$.

If the cohomological dimension of $\Gamma$ is $d-k+1$, then $\partial\Gamma$ has topological dimension $d-k$.
Since we have embedded $\partial\Gamma-\{x_0\}$ into a $(d-k)$-dimensional manifold, this implies that
$\partial \Gamma$ has a manifold point, which guarantees that $\partial\Gamma=S^{d-k}$, see Kapovich-Benakli \cite[Thm. 4.4]{kapovich-benakli}.
Let
$$E=\bigcup_{x\in\partial \Gamma} S(\xi(x)^{(k)})\subset S(\mathbb R^d),$$
where $S(v)$ denotes the unit sphere in $V$. Notice that $E$ is a closed submanifold of $S(\mathbb R^d)$ of dimension
$(d-k)+(k-1)=d-1$, so $E=S^{d-1}$. The map $p:E\to\partial \Gamma$ is a fibre bundle with fibre $S^{k-1}$, so the classification of
sphere fibrations (see \cite{adams}) implies that $(d-1,k-1)$ is either $(3,1)$, $(7,3)$ or $(15,7)$.
\qed

\medskip

Sambarino conjectured that if $\Gamma\subset \mathsf{PSL}(d,\mathbb R)$ is Borel Anosov, i.e. $P_k$-Anosov for all $1\le k\le d-1$,
then $\Gamma$ has a finite index subgroup which is either a free group or the fundamental group of a closed surface. This surprising
conjecture was proved when $d\le 4$ by Canary-Tsouvalas \cite{canary-tsouvalas}, for $d=2 (\mathrm{mod}\ 4)$ 
by Tsouvalas \cite{tsouvalas-borel}, and for
$d=3,\ 4 \text{ or }5 (\mathrm{mod}\ 8)$ by Dey \cite{dey-borel}.

\section{Hyperconvex groups}

Pozzetti, Sambarino and Wienhard \cite{PSW1} studied the class of $(1,1,2)$-hyperconvex groups. A group
$\Gamma\subset\mathsf{PSL}(d,\mathbb K)$ is $(1,1,2)$-hyperconvex if
it is $P_{\{1,2,d-2,d-1\} }$-Anosov and
whenever $x,y,z\in\partial\Gamma$ are distinct, then
$$\xi(x)^{(1)}\oplus\xi(y)^{(1)}\oplus \xi(z)^{(d-2)}=\mathbb K^d$$
where $\xi(x)^k$ is the $k$-dimensional subspace in the partial flag $\xi(x)$.
Farre, Pozzetti and Viaggi \cite{FPV} called these representation $(d-1)$-hyperconvex.
The most basic examples of  $(1,1,2)$-hyperconvex groups in $\mathsf{PSL}(d,\mathbb C)$ are the images of a convex cocompact
group subgroup of $\mathsf{PSL}(2,\mathbb C)$ by an irreducible representation of $\mathsf{PSL}(2,\mathbb C)$ into
$\mathsf{PSL}(d,\mathbb C)$.

One of the key motivations for \cite{PSW1} was  that $\alpha_1$-Patterson-Sullivan measures of  $(1,1,2)$-hyperconvex groups behave more
like the conformal measures developed by  Patterson \cite{patterson} and Sulllivan \cite{sullivan-hd} for Kleinian groups than Patterson-Sullivan measures
for typical Anosov groups.
The following result of  Pozzetti, Sambarino  and Wienhard \cite{PSW1} generalizes the fact, due to Patterson \cite{patterson} and Sulllivan \cite{sullivan-hd}, 
that the Hausdorff dimension of the limit set
of a convex cocompact  Kleinian group $\Gamma$ agrees with $\delta^{\alpha_1}(\Gamma)$.

\begin{theorem}{\rm (Pozzetti-Sambarino-Wienhard \cite{PSW1})} 
\label{haus dim}
If $\Gamma\subset\mathsf{PSL}(d,\mathbb K)$ is $(1,1,2)$-hyperconvex,
then the Hausdorff dimension of $\Lambda_1(\Gamma)$  is equal to $\delta^{\alpha_1}(\Gamma)$ and 
$\delta^{\alpha_1}(\Gamma)\le \mathrm{dim}(\mathbb K)$.
\end{theorem}

Pozzetti-Sambarino-Wienhard \cite[Prop. 4.1]{PSW1} and Glorieux-Monclair-Tholozan \cite[Thm 1.1]{GMT} showed that
if $\Gamma\subset\mathsf{PSL}(d,\mathbb K)$ is $P_{1,d-1}$-Anosov, then its limit set $\Lambda_1(\Gamma)$  has Hausdorff dimension at most $\delta^{\alpha_1}(\Gamma)$ but
that the bound is not always sharp.

\medskip

One of my favorite examples of the use of technology from Kleinian groups in higher rank is the following structural result of
Farre, Pozzetti and Viaggi.

\begin{theorem}{\rm (Farre-Pozzetti-Viaggi \cite[Thm. A]{FPV}) }
\label{hyperconvex char}
If $\Gamma\subset\mathsf{PSL}(d,\mathbb C)$ is $(1,1,2)$-hyperconvex, then $\Gamma$
is virtually isomorphic to a convex cocompact subgroup of $\mathsf{PSL}(2,\mathbb C)$.
\end{theorem}

We provide a brief outline of their beautiful  proof.
Given $x\in \partial \Gamma$, consider the projection map $\pi_x:\partial\Gamma\to \mathbb P(\mathbb C^d/\xi(x)^{(d-2)})\cong\mathbb{CP}^1$ 
(first studied in \cite{PSW1})
given by 
$$\pi_x(z)=\xi(z)^{(1)}\quad\text{if } z\ne x \quad\text{and}\quad\pi_x(x)=\xi(x)^{(d-1)}.$$
They show that this map is well-defined, continuous and injective. (The main difficulty here is establishing continuity at $x$.)

They then show, by a compactness argument, that there exists $K>1$ so that if $x,y\in \partial \Gamma$, then the homeomorphism
$$f_{x,y}=\pi_y\circ\pi_x^{-1}:\pi_x(\partial\Gamma)\to \pi_y(\partial \Gamma)$$
is $K$-quasi-M\"obius on $\pi_x(\partial\Gamma)$, i.e. if $z_1,z_2,z_3,z_4\in \pi_x(\partial\Gamma)$ and the  absolute value of their cross ratio $\big|[z_1,z_2,z_3,z_4]\big|=1$,
then 
$$\frac{1}{K}\le\Big| [f_{x,y}(z_1), f_{x,y}(z_2),f_{x,y}(z_3), f_{x,y}(z_4)]\Big|\le K.$$
If  $\partial\Gamma\cong S^2$, then $\xi_x(\partial\Gamma)=\mathbb P(\mathbb C^d/\xi(x)^{(d-2)})$ and one can conclude  that $f_{x,y}$ is 
$K$-quasiconformal  if it is orientation-preserving (since orientation-preserving $K$-quasi-M\"obius homeomorphisms of $\mathbb{CP}^1$ 
are $K$-quasiconformal.)

They then define an action  $\rho_x:\Gamma\to \mathrm{Homeo}(\pi_x(\partial \Gamma))$ given by
$$\rho_x(\gamma)=\pi_x\circ\gamma\circ (\pi_x)^{-1}\quad\text{for all}\ \ \gamma\in\Gamma$$
and check that
$$\rho_x(\gamma)=\hat\gamma|_{\gamma^{-}1(x),x}\circ f_{x,\gamma^{-1}(x)}$$
where $\hat\gamma_{\gamma^{-1}(x),x}:\mathbb P(\mathbb C^d/\xi(\gamma^{-1}(x))^{(d-2)})\to\mathbb P(\mathbb C^d/\xi(x)^{(d-2)})$  is
the map induced by $\gamma$.
Since $\hat\gamma_{\gamma^{-1}(x),x}$ is conformal, $\rho_x(\Gamma)$ is $K$-quasi-M\"obius for all $\gamma\in\Gamma$.

If $\partial\Gamma=S^2$, let
$$\Gamma_0=\{\gamma\in\Gamma:\rho_x(\gamma)\ \text{ is orientation-preserving }\}$$
and notice that $\Gamma_0$ has index at most two in $\Gamma$. So $\rho_x(\Gamma_0)$ is a uniformly
$K$-quasiconfomal action on $\mathbb P(\mathbb C^d/\xi(x)^{(d-2)})$, i.e. every $\rho_x(\gamma)$ is $K$-quasiconformal.
A theorem of Sullivan \cite[Theorem VII]{sullivan-ergodic}
then implies that $\rho_x(\Gamma_0)$ is quasiconformally conjugate to an action $\hat \Gamma$ of $\rho_x(\Gamma_0)$
on $\mathbb P(\mathbb C^d/\xi(x)^{(d-2)})$ by M\"obius transformations. Since $\Gamma$ acts on $\partial \Gamma$ as a uniform convergence group,
$\hat\Gamma$ acts on $\mathbb P(\mathbb C^d/\xi(x)^{(d-2)})$ as a uniform convergence group, so $\hat\Gamma$ is a uniform lattice
in $\mathsf{PSL}(2,\mathbb C)$.
If $K$ is the kernel of $\rho_x$, then $K$ is finite and $\hat\Gamma\cong \Gamma_0/K$. Therefore,
if $\partial\Gamma=S^2$, then $\Gamma$ is virtually isomorphic to a uniform lattice in $\mathsf{PSL}(2,\mathbb C)$.

In order to handle the general case, they use work of Haissinsky \cite{haissinsky} which implies  that a hyperbolic group $\Gamma$ 
with planar Gromov boundary is virtually isomorphic to a convex cocompact Kleinian group if every quasiconvex subgroup $H$  of $\Gamma$ 
whose Gromov boundary is homeomorphic to a Sierpinski carpet is  virtually isomorphic to a convex cocompact Kleinian group. 
In the case that $\partial H$
is a Sierpinski carpet,  they use results of  Ahlfors \cite{ahlfors-crit}, Bonk \cite{bonk}  and Markovic \cite{markovic}
to extend the uniformly $K$-quasi-M\"obius action $\rho_x(H)$ on 
$\pi_x(\partial H)$ to a uniformly $K'$-quasiconformal action on $\mathbb P(\mathbb C^d/\xi(x)^{(d-2)})$. 
(I am glossing over a delicate argument here.)
Sullivan's theorem again implies that
$\rho_x(H)$ is virtually isomorphic to a convex cocompact Kleinian group. This completes the (sketch of the) proof of Theorem \ref{hyperconvex char}.

\medskip

One may use the theory developed in Canary-Zhang-Zimmer \cite{CZZ2,CZZ3} to show that the $\alpha_1$-Patterson-Sullivan measure 
$\mu$ of a $(1,1,2)$-hyperconvex group is $\delta^{\alpha_1}$-Ahlfors regular, i.e. there exists $\epsilon_0,C>0$ so that
if $0<r<\epsilon_0$ and $F\in \Lambda_1(\Gamma)$, then 
$$\frac{1}{C}r^{\delta^{\alpha_1}(\Gamma)}\le \mu(B(F,r))\le  Cr^{\delta^{\alpha_1}(\Gamma)}$$
where $B(F,r)$ is the ball of radius $r$ about $F$ (in some fixed Riemannian metric on $\mathbb{P}(\mathbb K^d))$.
Since this implies that $\Lambda_1(\Gamma)$ has finite, non-zero $\delta^{\alpha_1}(\Gamma)$-dimensional Hausdorff measure,
one may regard this result as a strengthening of Theorem \ref{haus dim}.

\begin{theorem}{\rm (Canary-Zhang-Zimmer \cite{CZZ5})}
\label{ahlfors regular}
If $\Gamma\subset\mathsf{PSL}(d,\mathbb C)$ is $(1,1,2)$-hyperconvex, then the $\alpha_1$-Patterson-Sullivan measure $\mu$
of $\Gamma$ on $\Lambda_1(\Gamma)$ is $\delta^{\alpha_1}(\Gamma)$-Ahlfors regular. Therefore, $\Lambda_1(\Gamma)$ has
finite, non-zero $\delta^{\alpha_1}(\Gamma)$-dimensional Hausdorff measure.
\end{theorem}

Farre, Pozzetti and Viaggi   \cite[Thm 3.1]{FPV} proved that if $\partial\Gamma\ne S^2$, then $\pi_x(\partial \Gamma)$ has measure zero.
Their proof is a foliated version of Ahlfors'  \cite{ahlfors} original beautiful argument that the limit set of a geometrically
finite group has measure zero if it is not all of $\mathbb{CP}^1$.
It is easy to see that $\pi_x\circ \left(\xi^{(1)}\right)^{-1}:\Lambda_1(\Gamma)\to \pi_x(\partial\Gamma)$
is bilipschitz on any compact subset of $\Lambda_1(\Gamma)-\{\xi(x)^{(1)}\}$ (see the proof of \cite[Prop. 5.3]{FPV}). So,
if $\partial\Gamma\ne S^2$, then $\Lambda_1(\Gamma)$ has 2-dimensional Hausdorff measure zero. So, Theorem \ref{ahlfors regular}
implies that if $\delta^{\alpha_1}(\Gamma)=2$, then $\Gamma$ is virtually isomorphic to a lattice in $\mathsf{PSL}(2,\mathbb C)$.

\begin{theorem}{\rm (Canary-Zhang-Zimmer \cite{CZZ5})}
\label{2 implles lattice}
If $\Gamma\subset\mathsf{PSL}(d,\mathbb C)$ is $(1,1,2)$-hyperconvex, then
\hbox{$\delta^{\alpha_1}(\Gamma)=2$} if and only if $\Gamma$ is virtually isomorphic to a uniform lattice in
$\mathsf{PSL}(2,\mathbb C)$.
\end{theorem}

A result of Pozzetti, Sambarino and Wienhard \cite[Thm 7.1]{PSW1} (see also Zhang-Zimmer \cite{zhang-zimmer}) implies that if $\Gamma$ is
$(1,1,2)$-hyperconvex and  $\partial\Gamma=S^2$, then $\Lambda_1(\Gamma)$ is $C^1$. Since the tangent space  to $\Lambda_1(\Gamma)$
is a complex subspace at every point,  $\Lambda_1(\Gamma)$  is a complex submanifold, and hence algebraic, by a theorem of Chow \cite{chow}.
Therefore, the Zariski closure $Z$ of $\Gamma$ preserves $\Lambda_1(\Gamma)$ and acts as a group of biholomorphisms.
If $\Gamma$ is strongly irreducible, then $\Lambda_1(\Gamma)$ contains a projective frame (see \cite[Lemma 2.17]{BCLS}),
so the action of an element of $Z$ on $\mathbb P(\mathbb C^d)$ is determined by its action on $\Lambda_1(\Gamma)$.
Hence, $Z$ is an irreducible copy of $\mathsf{PSL}(2,\mathbb C)$ and we obtain the following rigidity theorem.

\begin{theorem}{\rm (Canary-Zhang-Zimmer \cite{CZZ5})}
\label{complex rigidity}
If $\Gamma\subset\mathsf{PSL}(d,\mathbb C)$ is $(1,1,2)$-hyperconvex and strongly irreducible and $\delta^{\alpha_1}(\Gamma)=2$, then
$\Gamma$  is the image of a  uniform lattice  in $\mathsf{PSL}(2,\mathbb C)$
by an irreducible representation of $\mathsf{PSL}(2,\mathbb C)$ into $\mathsf{PSL}(d,\mathbb C)$. 
\end{theorem}

Menal-Ferrer and Porti \cite{MP,porti} previously showed that if $\Gamma$ is a uniform lattice in $\mathsf{PSL}(2,\mathbb C)$ and 
$\tau_d:\mathsf{PSL}(2,\mathbb C)\to\mathsf{PSL}(d,\mathbb C)$ is an irreducible representation, then every small deformation
of $\tau_d(\Gamma)$ is conjugate to $\tau_d(\Gamma)$. Since every small deformation of a $(1,1,2)$-hyperconvex group
is $(1,1,2)$-hyperconvex, one may view Theorem \ref{complex rigidity} as a global version of the local rigidity
theorem of Menal-Ferrer and Porti.

\medskip\noindent
{\bf Historical remarks:} 
The results in this section also hold for $k$-hyperconvex groups, see \cite{FPV} and  \cite{CZZ5}.  A group
$\Gamma\subset\mathsf{PSL}(d,\mathbb K)$ is $k$-hyperconvex if
it is $P_{\{k-1,k,k+1,d-k-1,d-k,d-k+1\} }$-Anosov and
whenever $x,y,z\in\partial\Gamma$ are distinct, then
$$\Big( (\xi(x)^{(d-k)}\cap \xi(z)^{(k+1)})\oplus\xi(z)^{(k-1)}\Big)
\cap \Big( (\xi(y)^{(d-k)}\cap \xi(z)^{(k+1)})\oplus\xi(z)^{(k-1)}\Big) = \xi(z)^{(k-1)}.$$
Notice that a group is $(d-1)$-hyperconvex if and only if it is $(1,1,2)$-hyperconvex.
With the exception of Theorem \ref{complex rigidity}, one may recover these more general results quickly from
the results for $(1,1,2)$-hyperconvex groups by a careful study of the exterior power representation. The proof of
Theorem \ref{complex rigidity} requires additional Lie-theoretic arguments in the general case. We originally had a more complicated dynamical
proof of Theorem \ref{complex rigidity} but Andres Sambarino pointed out that we could use Chow's theorem to simplify the proof.

One can define a notion of $(1,1,2)$-hypertransversality for finitely generated transverse groups and prove that
if $\Gamma$ is $(1,1,2)$-hypertransverse, then $\delta^{\alpha_1}(\Gamma)$ is the Hausdorff dimension of the conical limit points
in $\Lambda_1(\Gamma)$, see \cite[Thm. 8.1]{CZZ2}. 
The proof contains ideas from a proof of analogous result for all finitely generated Kleinian groups
by Bishop and Jones \cite{bishop-jones} with ideas from the work of Pozzetti-Sambarino-Wienhard \cite{PSW1}.

Suppose that $\Gamma\subset\mathsf{PSL}(d,\mathbb K)$ is $P_\theta$-Anosov, $\phi\in\mathfrak{a}_\theta^*$ and $\delta^\phi(\Gamma)<+\infty$.
If $\phi$ is symmetric with respect to the obvious involution on $\mathfrak{a}_\theta^*$ (which takes $\omega_k$ to $\omega_{d-k}$),
then Dey and Kapovich \cite{DK-PS}  construct a Gromov pre-metric on $\mathcal F_\theta$ so that 
$\Lambda_\theta(\Gamma)$ has Hausdorff dimension $\delta^\phi(\Gamma)$. Dey, Kim and Oh \cite{DKO} proved that in
this same setting that the $\phi$-Patterson-Sullivan measure is $\delta^\phi(\Gamma)$-Ahlfors regular with respect to the pre-metric. 

In a sequel paper \cite{FPV2}, Farre, Pozzetti and Viaggi produce an analogue of Bers' Simultaneous Uniformization Theorem \cite{bers-simul}
for spaces of marked $k$-hyperconvex groups isomorphic to the fundamental group of a closed surface.  They also
show that if  $\Gamma\subset\mathsf{PSL}(d,\mathbb C)$ is fully hypererconvex, i.e. $k$-hyperconvex for all $k$,
and isomorphic to the fundamental group of a closed surface,
then its full limit set $\Lambda_{\{1,\ldots,d-1\} }$ has Hausdorff dimension 1 if and only if it is conjugate into $\mathsf{PSL}(d,\mathbb R)$.
One may view this last result as a generalization of Bowen's famous rigidity theorem \cite{bowen},
which asserts that the limit set of a quasifuchsian group
has Hausdorff dimension 1 if and only if it is Fuchsian.

\section{Combination theorems and their consequences}

Combination theorems arose in the study of Kleinian groups as a way of building new Kleinian groups from old.
The first general combination theorem was stated by Klein \cite{klein} in 1883. We recall that if $\Gamma$ is a discrete subgroup
of $\mathsf{PSL}(2,\mathbb C)$, then $\Gamma$ acts properly discontinuously on the complement $\Omega(\Gamma)$, in $\mathbb{CP}^1$, 
of its limit set $\Lambda(\Gamma)$. The set $\Omega(\Gamma)$ is called the domain of discontinuity and may be empty.

\medskip\noindent
{\bf Klein's combination theorem:} {\em Suppose that $\Gamma_1$ and  $\Gamma_2$ are  discrete subgroups of
$\mathsf{PSL}(2,\mathbb C)$ with non-empty domains of discontinuity. If $D_i$ is a fundamental domain for the action of
$\Gamma_i$ on $\Omega(\Gamma_i)$ and the closure of the exterior of $D_i$ is contained in $D_{i+1}$ (where we interpret
the indices modulo 2), then the group $\Gamma=\langle\Gamma_1,\Gamma_2\rangle$ generated by $\Gamma_1$ and $\Gamma_2$
is discrete and equal to $\Gamma_1*\Gamma_2$. Moreover, $D_1\cap D_2$ is a fundamental domain for the action of $\Gamma$
on $\Omega(\Gamma)$.}

\medskip

The proof of this theorem has now been immortalized as the ping-pong lemma. Notice that the analogous theorem where one
considers fundamental domains for the action on $\mathbb H^3$ also holds (with the same argument).

When one attempts to generalize this to the setting of Anosov groups, technical difficulties arise due to the fact
that $\Gamma$ does not act properly discontinuously on the complement of the limit set. However, Dey and Kapovich
establish the following analogue.

\begin{theorem} {\rm (Dey-Kapovich \cite{DK-FP})}
\label{klein anosov}
Suppose that $\Gamma_1$ and $\Gamma_2$ are $P_\theta$-Anosov  subgroups of $\mathsf{PSL}(d,\mathbb K)$,
that $A_i$ (for $i=1,2$) is a compact subset of $\mathcal F_\theta$ with non-empty interior and that every flag in $A_i$ is
transverse to every flag in $A_{i+1}$. If $\gamma_i(A_{i+1})\subset A_i$ for all $\gamma_i\in\Gamma_i$ (for  $i=1,2$),
then  $\langle\Gamma_1,\Gamma_2\rangle$
is $P_\theta$-Anosov and equal to $\Gamma_1*\Gamma_2$. 
\end{theorem}

In this formulation, one should think of $A_i$ as playing the role of $\mathbb{CP}^1-\mathrm{int}(D_i)$ in the original Klein
combination theorem. 

As a corollary one can show that any free product of Anosov groups is isomorphic to an Anosov group.

\begin{corollary} {\rm (Douba-Tsouvalas \cite[Thm. 3]{douba-tsouvalas}, Danciger-Gu\'eritaud-Kassel \cite[Cor 1.26]{DGK-combo})}
\label{free product of Anosov}
Suppose that $\Gamma_i$ is a $P_{\theta_i}$-Anosov  subgroups of $\mathsf{PSL}(d_i,\mathbb K)$, then there exists
a $P_{1,d-1}$-Anosov subgroup $\Gamma$ of $\mathsf{PSL}(d,\mathbb K)$ so that $\Gamma$ is isomorphic to $\Gamma_1*\Gamma_2$.
\end{corollary}

Danciger, Gu\'eritaud and Kassel \cite{DGK-combo} derive their proof of Corollary \ref{free product of Anosov} as
a consequence of their more general combination theorems for projectively convex cocompact, but not necessarily Anosov
or even hyperbolic, subgroups of $\mathsf{PSL}(d,\mathbb R)$. Douba and Tsouvalas \cite[Thm 3.1]{douba-tsouvalas} used
an early version of Theorem \ref{klein anosov} due to Dey-Kapovich-Leeb \cite{DKL} in their proof.

\medskip

Douba and Tsouvalas use Corollary \ref{free product of Anosov} in their proof that there exist Anosov groups not admitting
discrete faithful representations into any free product of rank one Lie groups.

\begin{theorem}{\rm (Douba-Tsouvalas \cite[Thm. 1/2]{douba-tsouvalas})}
Suppose that  $\Gamma_1$ is a uniform lattice in $F_4^{(-20)}$  and $\Gamma_2$ is a uniform lattice in $F_4^{(-20)}$  
or $\mathsf{Sp}(n,1)$ with $n\ge 51$. Then $\Gamma$ is isomorphic to a $P_{1,d-1}$-Anosov subgroup of $\mathsf{PSL}(d,\mathbb R)$,
for some $d$, but does not admit a discrete faithful representation into any finite product of rank one Lie groups.
\end{theorem}

\medskip

Maskit generalized Klein's combination theorem to allow for amalgamation along cyclic groups, free groups and surface groups.
We say that a closed subset $A\subset\mathbb{CP}^1$ is precisely invariant for a subgroup $H$ of  a Kleinian group $\Gamma$
if $h(A)=A$ for all $h\in H$ and $\gamma(A)\cap A=\emptyset$ for all $\gamma\in \Gamma-H$. We state Maskit's theorem
in the simpler setting of convex cocompact groups.

\begin{theorem}{\rm (Maskit \cite{maskit})}
\label{maskit combo}
Suppose that $\Gamma_1$ and  $\Gamma_2$ are  convex cocompact subgroups of
$\mathsf{PSL}(2,\mathbb C)$ and that $H=\Gamma_1\cap\Gamma_2$ is a  convex cocompact subgroup which has infinite
index in both $\Gamma_1$ and $\Gamma_2$. Suppose that the limit set $\Lambda(H)$ is contained in
a Jordan curve $J$, and let $A_1$ and $A_2$ be the components of $\mathbb CP^1-J$.
If  $A_i$ is precisely invariant for $H$ in $\Gamma_i$ (for $i=1,2$), then
$\Gamma=\langle\Gamma_1,\Gamma_2\rangle$ is discrete and isomorphic to $\Gamma_1*_H\Gamma_2$.
\end{theorem}

Maskit's theorem played  a crucial role in Thurston's proof of his hyperbolization theorem, see Morgan \cite{morgan}.

Dey and Kapovich also obtained an analogue of Maskit's theorem. Moreover, Danciger, Gu\'eritaud and Kassel \cite{DGK-combo}
establish a version for projectively convex cocompact groups.

\begin{theorem} {\rm (Dey-Kapovich \cite{DK-KMA})}
\label{Anosov amalgam}
Suppose that $\Gamma_1$ and $\Gamma_2$ are $P_\theta$-Anosov  subgroups of $\mathsf{PSL}(d,\mathbb K)$,
$H=\Gamma_1\cap\Gamma_2$ is quasiconvex in either $\Gamma_1$ or $\Gamma_2$. If $A_1$ and $A_2$ are compact subsets of
$\mathcal F_\theta$ so that
\begin{enumerate}
\item
Every flag in the interior of  $A_1$ is transverse to every flag in the interior of $A_2$,
\item
$A_i$ is precisely invariant for $H$ in $\Gamma_i$ and
\item
Every flag in $A_i$ is transverse to every flag in $\Lambda_\theta(\Gamma_{i+1})-\Lambda_\theta (H)$,
\end{enumerate}
then $\langle\Gamma_1,\Gamma_2\rangle$
is $P_\theta$-Anosov and equal to $\Gamma_1*_H\Gamma_2$. 
\end{theorem}

Dey and Tsouvalas used Theorem \ref{Anosov amalgam} and a separability result due to Tsouvalas 
\cite{tsouvalas-sep}  to produce many new Anosov groups by amalgamating along
cyclic subgroups. The following example is especially relevant.

\begin{theorem}{\rm (Dey-Tsouvalas \cite{dey-tsouvalas})}
If $\Gamma$ is a uniform lattice in $\mathsf{Sp}(n,1)$ and $H$ is an infinite abelian subgroup of $\Gamma$, then there
exists a finite index subgroup $\Gamma'$ of $\Gamma$ containing $H$, so that $\Gamma'*_H\Gamma'$ admits a
$P_{1,d-1}$-Anosov representation into $\mathsf{SL}(d,\mathbb C)$ for some $d$.
\end{theorem}

Combining this with earlier work of Tholozan and Tsouvalas \cite{tholozan-tsouvalas} produces many examples of
one-ended hyperbolic groups which are isomorphic to Anosov groups, but do not admit discrete, faithful representations
into any rank one Lie group.

\begin{corollary}{\rm (Dey-Tsouvalas \cite{dey-tsouvalas})}
There exist one-ended Anosov groups which do not admit any discrete, faithful representation into
a rank one Lie group.
\end{corollary}

\medskip\noindent
{\bf Historical remarks:}
Traaseth and Weisman \cite{traaseth-weisman} proved combination theorems for geometrically finite convergence group actions. Tsouvalas
and Weisman  \cite{tsouvalas-weisman} proved combination theorems for quasi-isometrically embedded groups. As a consequence
they see that if $\Gamma_1$ and $\Gamma_2$ are quasi-isometrically embedded subgroups of $\mathsf{PSL}(d_1,\mathbb K)$ and
$\mathsf{PSL}(d_2,\mathbb K)$, then $\Gamma_1*\Gamma_2$ admits a quasi-isometric embedding into some $\mathsf{PSL}(d,\mathbb K)$.
Both Dey-Kapovich \cite{DK-KMA} and Danciger-Gu\'eritaud-Kassel \cite{DGK-combo} also establish combination
theorems for HNN extensions in the spirit of Theorem \ref{maskit combo}.

Recall that a domain $\Omega$ in $\mathbb P(\mathbb R^d)$ is properly convex if it is a convex bounded subset of some affine
chart for $\mathbb P(\mathbb R^d)$. If a  subgroup $\Gamma\subset\mathsf{PSL}(d,\mathbb R)$ preserves  and acts properly
discontinuously on a  properly convex
domain, then its full orbital limit set $\Lambda_\Omega(\Gamma)\subset\partial \Omega$ is the set of accumulation points of any orbit, i.e.
$z\in \Lambda_\Omega(\Gamma)$ if and only if there exists $x\in\Omega$ and $\{\gamma_n\}\subset\Gamma$ so that $\gamma_n(x)\to z$.
We say that $\Gamma$ is (projectively) convex cocompact if $\Gamma$ acts cocompactly on the convex hull, 
in $\Omega$, of $\Lambda_\Omega(\Gamma)$. 
Convex cocompact Kleinian groups are also convex cocompact in this definition where we take $\Omega$ to be the round ball in
$\mathbb P(\mathbb R^4)$ preserved by $\mathsf{PO}(3,1)$.  More generally, any Anosov subgroup of $\mathsf{PSL}(d,\mathbb R)$
which preserves a properly convex domain is projectively convex cocompact, see Danciger-Gu\'eritaud-Kassel \cite[Thm 1.4]{DGK-CC}
or Zimmer \cite[Thm 1.27]{zimmer-cc}.
However, projectively convex cocompact groups need not be Anosov or
even word hyperbolic.

\section{Proper affine actions}

In this section we will discuss  the affine group $\mathrm{Aff}(V)$ of a finite-dimensional $\mathbb K$-vector space $V$.
The affine group  $\mathrm{Aff}(V)$ may be identified with the semidirect product of $\mathsf{GL}(V)$ and $V$,
where the action is given by
$$(A,\vec v)(\vec w)=A(\vec w)+\vec v.$$
The affine group is not semi-simple so does not fit into the framework of the previous sections, but its development has been
heavily influenced by ideas from Fuchsian and Kleinian groups. For a fuller treatment of this subject we recommend the
survey paper by Danciger, Drumm, Goldman and Smilga \cite{DDGS}.

Auslander \cite{auslander} conjectured that if $\Gamma\subset \mathrm{Aff}(\mathbb R^d)$ acts properly discontinuously and cocompactly 
on $\mathbb R^d$ by affine transformations, then $\Gamma$ is virtually solvable. The conjecture remains open, 
but has been proven in all dimensions  up to six, see  Fried-Goldman \cite{fried-goldman}, Tomanov \cite{tomanov} and Abels-Margulis-Soifer \cite{AMS}.
Milnor \cite{milnor} asked whether something similar might be true for actions which are not cocompact, in analogy with the
Bieberbach theorems for groups of Euclidean isometries.

Margulis \cite{Mar1,Mar2} produced the first examples of proper affine actions by non-abelian free groups on $\mathbb R^3$, which
are now called Margulis space-times. His examples may be viewed as arising from a one-parameter family of  convex cocompact representations
$\{\rho_t:F_n\to \mathsf{SL}(2,\mathbb R)\}$ where we identify $\mathfrak{sl}(2,\mathbb R)$ with $\mathbb R^3$ and $\gamma\in F_n$
acts by
 $$\gamma\to\left(\mathrm Ad (\rho_0(\gamma)),\left(\frac{d}{dt}\Big|_{t=0}\rho_t(\gamma)\right)\rho_0(\gamma)^{-1}\right)$$
 where $\mathrm{Ad}:\mathsf{SL}(2,\mathbb R)\to\mathfrak{sl}(2,\mathbb R)$ is the adjoint representation.
 The crucial geometric  feature of the deformation is that if $\ell(\rho_t(\gamma))$ denotes the translation length of $\rho_t(\gamma)$, then
 there exists $c>0$ so that 
 $$\frac{d}{dt}\Big|_{t=0} \ell(\rho_t(\gamma))\le -c\ell(\rho_0(\gamma))$$
 for all $\gamma\in F_n$.

Drumm \cite{drumm} and Drumm-Goldman \cite{DG-crooked} introduced a geometric viewpoint
on Margulis' construction, produced large classes of new examples and exhibited fundamental domains for their examples.
Danciger, Gu\'eritaud and Kassel \cite{DGK1,DGK2}
gave a complete classification of Margulis space-times with convex cocompact linear part and showed that their quotients are all
homeomorphic to the interior of a handlebody.

Goldman, Labourie and Margulis  \cite{GLM} gave 
an exact criterion for when affine actions of free groups on $\mathbb R^3$ are proper (see Ghosh-Trieb for $\mathbb R^{2n+1}$).
Smilga \cite{Smig} extended  Margulis' construction to any noncompact semisimple  Lie group $\mathsf{G}$. He
constructed a proper affine  action of a non-abelian free group on the Lie algebra $\mathfrak g$ whose  linear part is Zariski dense in
$\mathrm{Ad}(\mathsf{G})$.  Abels, Margulis and Soifer \cite{AMS-orthogonal} showed that if the linear part of a proper action
is a Zariski dense subgroup of $\mathsf{O}(p,q)$ with $p\ge q\ge 1$, then $(p,q)=(2n,2n-1)$.
Smilga  \cite{smilga-ortho}  described fundamental domains for some proper
affine actions of  free groups on $\mathbb R^{4n-1}$ whose linear part lies in  $\mathsf{SO}(2n,2n-1)$.

$\check{\mathrm{Z}}$ager Korenjak  \cite{neza} generalized the strip deformations of Danciger, Gu\'eritaud and Kassel
to produce further proper affine actions of  free groups on $\mathbb R^{4n-1}$ whose linear part lies in  $\mathsf{SO}(2n,2n-1)$.
Burelle and $\check{\mathrm{Z}}$ager Korenjak \cite{jp-neza} showed that the image of every positive representation of a free group
into $\mathsf{SO}(2n,2n-1)$ arises as the linear part of a proper action on $\mathbb R^{4n-1}$. They also exhibited fundamental domains
for actions arising from strip deformations which generalize Drumm's crooked plane description \cite{drumm}  from the $n=1$ case.

\medskip
 
It is natural to ask which other groups admit proper affine actions. In a major breakthrough, Danciger, Gu\'eritaud and Kassel \cite{DGK-aff} 
proved that many geometrically natural groups admit proper affine actions, including surface groups, all hyperbolic 3-manifold groups,
and groups of arbitrarily large cohomological dimension.

\begin{theorem}{\rm (Danciger-Gu\'eritaud-Kassel \cite{DGK-aff})} 
Any right-angled Coxeter group on $k$ generators admits a proper affine action on $\mathbb R^{k(k-1)/2}$.
\end{theorem}

In this setting, they construct families $\{\rho_t\}$ of representations  of the right-angled Coxeter group into $\mathsf{O}(p,q+1)$ where $p+q+1=k$
and obtain actions on $\mathfrak{o}(p,q+1)$ where the linear part is the (image under the adjoint representation) of $\rho_0$ and
the translational part is given by the derivative of the deformation. They develop a contraction property which guarantees that the resulting
affine action is proper and verify their criterion holds for their construction.

Danciger, Gu\'eritaud and Kassel \cite[Prop. 1.6]{DGK-aff} also show that every hyperbolic surface group admits a proper affine action 
on $\mathbb R^6\cong\mathfrak{sl}(2,\mathbb C)$. For comparison,  Mess  \cite{Mess} had previously shown that a closed surface group 
cannot have a proper affine action on $\mathbb R^3$. In the surface case, they begin with a  subgroup of $\mathrm{Isom}(\mathbb H^2)$
which is the reflection group of a right-angled $(2g+2)$-gon $P$ in $\mathbb H^2$, 
which contains the fundamental group of a closed surface of genus $g$ as an index 4 subgroup.
They extend $P$ to a right-angled polyhedron $\hat P$ in $\mathbb H^3$ which also has $2g+2$ faces. They explicitly describe a
smooth one-parameter family of  deformations of $\hat P$ within $\mathbb H^3$ which remain right-angled and hence generate reflection
groups in $\mathbb H^3$. They show that  in a small neighborhood of the Fuchsian group (but not at the Fuchsian group itself) the derivatives of
the elements of the group have the contraction property guaranteeing that they give rise to  proper affine actions.

\medskip

Bridgeman, Canary and Sambarino \cite{BCS} recently used the deformation theory of quasifuchsian groups to produce
an open subset $U$ of the space of quasifuchsian groups such that (the image under the adjoint) of every quasifuchsian group in $U$
is the linear part of a proper affine action on $\mathfrak{sl}2,\mathbb C)\cong\mathbb R^6$. We recall that the space $QF(S)$ of
marked quasi-fuchsian groups isomorphic to $\pi_1(S)$ is the space of (conjugacy classes of) convex compact representations
of $\pi_1(S)$ into $\mathsf{PSL}(2,\mathbb C)$. Bers \cite{bers-simul} showed that $QF(S)$ is naturally identified with 
$\mathcal T(S)\times\mathcal T(\bar S)$ where $\mathcal T(S)$ is the Teichm\"uller space of marked conformal
structures on $S$ and $\bar S$ is $S$ with the opposite orientation. The space of Fuchsian representations (i.e. representations
conjugate into $\mathsf{PSL}(2,\mathbb R)$) manifests as the diagonal in this parametrization.

If $\rho\in QF(S)$, then its limit set $\Lambda(\rho)=\Lambda(\rho(\Gamma))\subset \partial\mathbb H^3$ is a Jordan curve.
Let $CH(\rho)$ be the convex hull in $\mathbb H^3$ of $\Lambda(\rho)$ and $C(\rho)=CH(\rho)/\rho(\pi_1(S))$
be the convex core of $N_\rho=\mathbb H^3/\rho(\pi_1(S))$.  Then $C(\rho)$ is homeomorphic to $S\times [0,1]$ unless
$\rho$ is Fuchsian.

Thurston \cite{thurston-notes} (see also Epstein-Marden \cite{epstein-marden}) showed that each component of 
the boundary of $C(\rho)$ is a hyperbolic surface in its intrinsic metric. Moreover, Thurston showed that each
component of the boundary of $C(\rho)$ is totally geodesic in the complement of a lamination called the bending
lamination. Moreover, the bending lamination inherits a transverse measure which measures the total bending along
an arc in the surface. If $\rho$ is not Fuchsian one obtains a pair $(\beta_+,\beta_-)$ of measured laminations. Bonahon
and Otal \cite{bonahon-otal} characterized exactly which pairs of laminations can arise. Dular and Schlenker \cite{dular-schlenker} recently showed that the
pair $(\beta_+,\beta_-)$ determines $\rho$.

If a simple closed curve $\gamma$ separates $S$ into two surfaces $S_1$ and $S_2$, then
$$\pi_1(S)=\pi_1(S_1)*_{\langle\gamma\rangle}\pi_1(S_2).$$
Given $\rho\in QF(S)$ one can define a bending deformation of $\rho$ along $\gamma$. Let $A_{\rho(\gamma)}(\theta)$
be the M\"obius transformation which rotates by an angle $\theta$ in the axis of $\rho(\gamma)$. We then
define $\rho_\theta:\pi_1(S)\to\mathsf{PSL}(2,\mathbb C)$ by $\rho_\theta=\rho$ on $\pi_1(S_1)$ and
$\rho_\theta=A_{\rho(\gamma)}(\theta)\rho A_{\rho(\gamma)}(\theta)^{-1}$ on $\pi_1(S_2)$. Notice that $\rho_0=\rho$
and $\rho_\theta\in QF(S)$ for all small enough values of $\theta$.

Kourouniotis \cite{kouron} gave a formula for the derivative of the complex length of every element under the bending deformation along
a curve or multicurve.
One can, by taking limits, define a bending deformation of a quasifuchsian group along any measured lamination on the surface,
see Kourouniotis  \cite{kouron-cont} and Epstein-Marden \cite{epstein-marden}. One can generalize Kourouniotis' work
to give a formula for the derivative of the complex length of elements with respect to the bending deformation along a general measured lamination.
We use techniques developed by Bridgeman-Canary-Yarmola \cite{BCY} to develop a criterion
which guarantees that if $\{\rho_t\}$ is the bending deformation of $\rho$ along $\beta_+$ (or along $\beta_-$), then there exists $c>0$ so that
$$\frac{d}{dt}\Big|_{t=0} \ell(\rho_t(\gamma))\le -c i(\gamma,\beta_+)$$
where $\ell(\rho_t(\gamma))$ is the real translation length of $\rho_t(\gamma)$ and $ i(\gamma,\beta_+)$ is geometric
intersection number of $\gamma$ and $\beta_+$.

Let $\vec v_\pm\in T_\rho QF(S)$
be the tangent vectors given by the bending deformation along $\beta_\pm$. Set $\vec v=\vec v_+ +\vec v_-$ and let
$\{\rho_t\}$ be a deformation of $\rho$ in the direction $\vec v$. One can show that if both $\beta_+$ and $\beta_-$ satisfy
our earlier criterion then there exists $d>0$ so that 
\begin{equation}
\label{uniform shrink}
\frac{d}{dt}\Big|_{t=0} \ell(\rho_t(\gamma))\le -d\ell(\rho_0(\gamma))
\end{equation}
for all $\gamma\in\pi_1(S)$.

Results of Ghosh \cite{ghosh} and Kassel-Smilga \cite{kassel-smilga}, combined with recent work of Sambarino \cite{sambarino}
can then be used to guarantee that Equation \eqref{uniform shrink} implies that 
$$\gamma\to\left(\mathrm{Ad} (\rho(\gamma)),\left(\frac{d}{dt}\Big|_{t=0}\rho_t(\gamma)\right)\rho(\gamma)^{-1}\right)$$
gives a proper affine action of $\pi_1(S)$,
where $\mathrm{Ad}:\mathsf{SL}(2,\mathbb C)\to\mathfrak{sl}(2,\mathbb C)$ is the adjoint representation.
We have completed the outline of the  proof of  the following result.

\begin{theorem}{\rm (Bridgeman-Canary-Sambarino \cite{BCS})}
If $S$ is a closed surface of genus $g\ge 2$, then there exists an open neighborhood $U$ of the Fuchsian locus
in $QF(S)$ so that if $\rho\in U$ is not Fuchsian, then $\mathrm{Ad}(\rho)$ is the linear part of a proper affine
action of $\pi_1(S)$ on $\mathfrak{sl}(2,\mathbb C)\cong\mathbb R^6$.
\end{theorem}

\medskip

Our original motivation for studying bending deformations was the entropy functional on quasifuchsian space.
Here the entropy  of $\rho\in QF(S)$ is given by $\delta^{\alpha_1}(\rho(\pi_1(S))$.
Ruelle \cite{ruelle} showed that the entropy functional is analytic on $QF(S)$ and Bowen \cite{bowen} showed
that it achieves its global minimum of 1 only along the Fuchsian locus.  Bridgeman \cite{bridgeman-wp}
showed that  the Hessian of the entropy functional is positive definite on at least a half-dimensional subspace, so the entropy 
functional has no local maxima.

Work of Sambarino \cite{sambarino} shows that if $\{\rho_t\}$ is a smooth path  in $QF(S)$ and
$$\frac{d}{dt}\Big|_{t=0} \ell(\rho_t(\gamma))\le 0 \text{ for all } \gamma\in\pi_1(S) \text{ and } \frac{d}{dt}\Big|_{t=0} \ell(\rho_t(\alpha)< 0
\text{ for some }\alpha\in\pi_1(S),$$
then $\rho_0$ is not a critical point of the entropy functional.
Therefore, if $\rho\in QF(S)$ is not Fuchsian and  either bending lamination,   $\beta_+$  or $\beta_-$, satisfies our early criterion, 
then  we may bend along that bending lamination to show that $\rho$ is not a critical point of the energy functional.

\begin{theorem}{\rm (Bridgeman-Canary-Sambarino \cite{BCS})}
If $S$ is a closed surface of genus $g\ge 2$, then there exists an open neighborhood $U$ of the Fuchsian locus
in $QF(S)$ so that if $\rho\in U$ is not Fuchsian, then $\rho$ is not a critical point of the energy functional.
\end{theorem}

\section{Other advances}

In this section, we very briefly discuss a few other recent developments motivated by ideas from the world of Kleinian groups.

\medskip\noindent
{\bf Analogues of geometric finiteness:}  It is well-established that Anosov groups are the analogue in higher rank of convex cocompact
groups in rank one. It is then natural to ask for the higher rank analogue of geometrically finite subgroups of rank one Lie groups.
Here there are two competing notions.

Kapovich and Leeb \cite{KL} were the first to develop the theory of what is now known as relatively Anosov groups. This
theory was further developed by Zhu \cite{zhu-thesis} and Zhu-Zimmer \cite{ZZ1,ZZ2}. Given a relatively Anosov groups
one can show that if the critical exponent of a linear functional is finite, then its Poincar\'e series diverges at its critical exponent, see \cite{CZZ4}.
(The proof follows the outline of Dal'bo, Otal and Peign\'e \cite{DOP} but in higher rank substantial difficulties arise in showing
that the Poincar\'e series of any peripheral subgroup diverges at its critical exponent, and we must employ Hironaka's
theorem \cite{hironaka} on resolution of singularities  in the proof.)
One may further construct a BMS measure on the geodesic flow associated to the linear functional  which one can
show is finite and that the flow is mixing, see Blayac-Canary-Zhu-Zimmer \cite{BCZZ2} and  Kim-Oh \cite{KO-RA}. One may
thus derive counting theorems, see \cite{BCZZ2}. Counting and mixing results were previously established for images of relatively Anosov
representations of geometrically finite Fuchsian groups, see Bray-Canary-Kao-Martone \cite{BCKM}, for rank one properly convex
projective structures, see Blayac-Zhu \cite{BZ}, and for properly (but not strictly) convex  closed projective 3-manifolds,
see Bray \cite{bray1,bray2}.

The main disadvantage of the theory of relatively Anosov groups is that it does not cover many of the examples which one would
naturally like to call geometrically finite. The restrictions largely arise from the fact that relatively Anosov subgroups  of 
$\mathsf{PSL}(d,\mathbb K)$ are relatively
hyperbolic groups whose peripheral subgroups consist of weakly unipotent elements (i.e. elements all of whose eigenvalues have
modulus 1). This feature, along with the fact that there limit sets are transverse, allows for the dynamical analysis discussed
in the previous paragraph.

Weisman \cite{weisman-EGF} introduced the more general class of extended geometrically finite groups. His class of groups
still consists of relatively hyperbolic groups but their limit sets need not be transverse and the elements of their peripheral
subgroups need not be weakly unipotent. His class of groups includes all relatively Anosov groups,  all projectively convex cocompact groups which
are relatively hyperbolic and  holonomies of many convex projective manifolds with generalized cusps, see Weisman \cite{weisman-EGF2}.
(Islam and Zimmer \cite{IZ,IZ2} gave necessary and sufficient conditions for a projectively convex cocompact group to be relatively hyperbolic.)
Impressively, Weisman \cite{weisman-EGF} is able to prove a strong stability theorem in his setting. Extended geometrically finite groups
have not been as extensively studied as relatively Anosov groups, but seem likely to become an important organizing principle for
this larger class of groups.

\medskip\noindent
{\bf Limits of Anosov groups:} As the theory of Kleinian groups matured, it began to focus on the study of limits of geometrically finite
groups. This theory is in its infancy in the higher rank setting. Schwartz \cite{schwartz-pappus} constructed discrete subgroups
of $\mathsf{PSL}(3,\mathbb R)$ which are isomorphic to $\mathsf{PSL}(2,\mathbb Z)$. Barbot, Lee and Val\'erio \cite{BLV} showed
that these groups arise as algebraic limits of Anosov groups, and the limiting behavior was further analyzed by Schwartz \cite{schwartz-anosov}.
These groups are now known to be relatively Anosov. Lahn \cite{lahn1,lahn2} has extensively studied limiting behavior of families
of reducible representations.

Very recently, Bobb and Farre discovered representations of surface groups into $\mathsf{PSL}(4,\mathbb R)$ which are limits
of Anosov representations but not even extended geometrically finite (and are un-related to geometrically infinite Kleinian groups).
They associate an ``ending lamination'' to these representations which encodes the failure of the Anosov property. They
also analyze and draw pictures of the limit sets.
This discovery was an outgrowth of their work in \cite{bobb-farre} on the convex core of surface groups acting convex compactly on
$\mathbb P(\mathbb R^4)$. This discovery suggests many exciting new avenues for research.

It is known that every word hyperbolic Kleinian group is an algebraic  limit of convex cocompact Kleinian groups, see
Namazi-Souoto \cite{namazi-souto} and Ohshika \cite{ohshika-density}. Tsouvalas \cite{tsouvalas-density} decisively proved that the analogue
is not true in general in higher rank. He exhibited a hyperbolic group $\Gamma$ and an open set in 
$\mathrm{Hom}(\Gamma,\mathsf{PSL}(d,\mathbb K))$ (for both $\mathbb K=\mathbb R$ and $\mathbb K=\mathbb C$) consisting
of discrete, faithful representations which are quasi-isometric embeddings whose images are Zariski dense but not Anosov.

\medskip\noindent
{\bf Dehn filling:} Thurston \cite{thurston-notes} famously showed that all but finitely many Dehn fillings of a one-cusped
finite volume hyperbolic 3-manifolds admit hyperbolic structures. Moreover, as the Dehn surgery coefficients diverge to infinity
the resulting Dehn filled  hyperbolic manifolds converge to the original one-cusped hyperbolic 3-manifold. 
Choi, Lee and Marquis \cite{choi-lee-marquis} established a generalized Dehn filling theorem for certain convex projective manifolds
in dimensions 4 through 7. Schwartz \cite{schwartz-CR} and Acosta \cite{acosta} proved  Dehn filling theorems for certain classes of subgroups
of $\mathsf{SU}(2,1)$.

Weisman \cite{weisman-dehn} recently proved a surprisingly strong generalization of Thurston's theorem which includes Thurston's
results and the results of Choi-Lee-Marquis as special cases. (It appears likely that the results of Acosta and Schwartz are also
special cases.)  As a first application, in collaboration with Danciger,
he constructs exotic new representations of 3-manifolds groups into $\mathsf{PU}(3,1)$.

\medskip\noindent
{\bf Other topics:} Here I briefly mention a few topics that I ran out time, space or expertise to discuss properly.

\begin{enumerate}
\item
{\em Pleated surfaces:} Maloni, Martone, Mazzoli and Zhang \cite{MMMZ1,MMMZ2} have developed a theory of
$d$-pleated surfaces associated to representations of surface groups into $\mathsf{PSL}(d,\mathbb C)$. Mazzoli and Viaggi \cite{MV}
also develop a theory of pleated surfaces in $\mathbb H^{2,n}$ in their study of maximal representations into $\mathsf{SO}_0(2,n+1)$.
\item
{\em Surface subgroups:} Inspired by the solution of the surface subgroup problem for hyperbolic 3-manifolds by
Kahn and Markovic \cite{kahn-markovic}, Kahn, Labourie and Mozes \cite{KLM} showed that uniform lattices
in center-free complex semi-simple Lie groups contain surface subgroups.
\item
{\em Anti de Sitter 3-manifolds:} There is a long history of techniques from Teichm\"uller theory and the theory
of quasifuchsian manifolds influencing the study of anti de Sitter 3-manifolds, beginning with the seminal
work of Mess \cite{Mess}, see also \cite{mess-notes}. For a recent survey of related work see Bonsante-Seppi 
\cite{bonsante-seppi}.
\end{enumerate}


\begin{thebibliography}{100}
\begin{scriptsize}

\bibitem{AMS} H. Abels, G. Margulis and A. Soifer, ``Properly discontinuous groups of affine transformations with orthogonal linear part,''
{\em C. R. Acad. Sci. Paris  Math.} {\bf 324}(1997), no. 3, 253--258.

\bibitem{AMS-orthogonal} H. Abels, G. Margulis and A. Soifer, ``On the Zariski closure of the linear part of a properly discontinuous group of 
affine transformations,''
{\em J. Diff. Geom} {\bf 60}(2002),  315--344.

\bibitem{acosta} M Acosta, ``Spherical CR uniformization of Dehn surgeries of the Whitehead link complement,'' 
{\em  Geom. Top.} {\bf 23}(2019),  2593--2664. 

\bibitem{adams} J.F. Adams, ``On the non-existence of elements of Hopf invariant one,'' {\em Ann. Math.} {\bf 72}(1960), 20--104.

\bibitem{ahlfors-crit} L. Ahlfors,  ``Quasiconformal reflections,'' {\em  Acta Math.} {\bf 109}(1963), 291--301. 

\bibitem{ahlfors} L. Ahlfors, ``Fundamental polyhedrons and limit point sets of Kleinian groups,'' \emph{Proc. Nat. Acad. Sci. USA}
\textbf{55} (1966), 251--254.

\bibitem{albuquerque} P. Albuquerque, ``Patterson-Sullivan theory in higher rank symmetric spaces,'' {\em G.A.F.A.} {\bf  9}(1999), 1--28.

\bibitem{mess-notes} L. Andersson, T. Barbot, R. Benedetti, F. Bonsante, W. Goldman, F.  Labourie, K. Scannell and J.-M. Schlenker,
``Notes on a paper of Mess,'' {\em Geom. Ded.} {\bf 126}(2007), 47--70.


\bibitem{auslander} L. Auslander, ``The structure of complete locally affine manifolds,''
{\em Topology}  {\bf 3}(1964), 131--139.

\bibitem{BLV} T. Barbot, G.-S. Lee and V. Val\'erio, ``Pappus theorem, Schwartz representations and Anosov representations,''
{\em Ann Inst. Four.} {\bf 68}(2018), 2697--2741.

\bibitem{benoist-divisible1} Y. Benoist, ``Convexes divisibles I,'' in {\em Algebraic groups and arithmetic},
Tata Inst. Fund. Res. Stud. Math. {\bf 17}(2004), 339--374.

\bibitem{benoist-quint} Y. Benoist and J.F. Quint, {\em Random Walks on Reductive Groups}, Springer-Verlag, 2016.

\bibitem{bers-simul} L. Bers, ``Simultaneous uniformization,'' {\em Bull. A.M.S.} {\bf 66}(1960), 94--97. 

\bibitem{bestvina-mess} M. Bestvina and G. Mess, ``The boundary of a negatively curved group,''
{\em J.A.M.S.} {\bf 4}(1991), 469--481.

\bibitem{bishop-jones} C. Bishop and P. Jones, ``Hausdorff dimension and Kleinian groups,''  {\em Acta Math.} {\bf 179}(1997),  1---39.

\bibitem{BCZZ2} P.-L. Blayac, R. Canary, F. Zhu, and A. Zimmer, ``Counting, mixing and equidistribution for GPS 
systems with applications to relatively Anosov groups,'' preprint, arXiv:2404.09718.

\bibitem{BZ} P.-L. Blayac and F. Zhu, ``Ergodicity and equidistribution in Hilbert geometry,'' {\em J. Mod. Dyn.}  {\bf 19}(2023), 879--94.

\bibitem{bobb-farre} M. Bobb and J. Farre, ``Affine laminations and coaffine representations,'' preprint, arXiv:2404.14284.

\bibitem{BPS} J. Bochi, R. Potrie and A. Sambarino, ``Anosov representations and dominated splittings,'' {\em J.E.M.S.} {\bf 21}(2019), 3343--3414.

\bibitem{bonahon-otal} F. Bonahon and J. P. Otal, ``Laminations measur\'ees de plissage des vari\'et\'es hyperboliques de dimension 3,''
{\em Ann. of Math.} {\bf 160}(2004),  1013--1055.

\bibitem{bonk} M. Bonk, ``Uniformization of Sierpinski carpets in the plane,'' {\em Invent. Math.} {\bf186}(2011), 559--665.

\bibitem{bonsante-seppi} F. Bonsante and A. Seppi, ``Anti-de Sitter geometry and Teichm\"uller theory,'' in
{\em In the tradition of Thurston: Geometry and Topology}, Springer Verlag(2020), 543--563.


\bibitem{bowen} R. Bowen, ``Hausdorff dimension of quasi-circles",
{\em Pub. Math. I.H.E.S.} {\bf 50}(1979), 11--25.

\bibitem{bray1} H. Bray, ``Geodesic flow of nonstrictly convex Hilbert geometries,''
{\em Ann. Inst. Four.} {\bf 70}(2020), 1563--1593.

\bibitem{bray2} H. Bray, ``Ergodicity of Bowen-Margulis measure for the Benoist 3-manifolds,''
{\em J. Mod. Dyn.} {\bf 16}(2020)  305--329.

\bibitem{BCKM} H. Bray, R. Canary, L.-Y. Kao, and G. Martone, ``Counting, equidistribution and entropy gaps at infinity with applications to 
cusped Hitchin representations,''  {\em J. Reine Angew. Math.}  {\bf 791}(2022), 1--51.


 
\bibitem{bridgeman-wp} M. Bridgeman, ``Hausdorff dimension and the Weil-Petersson extension
to quasifuchsian space,''
{\em Geom. and Top.} {\bf 14}(2010),  799--831.

\bibitem{BCLS} M. Bridgeman, R. Canary, F. Labourie and A. Sambarino, ``The pressure metric for Anosov representations,'' {\em G.A.F.A.}, {\bf 25}(2015), 1089--1179.

\bibitem{BCS} M. Bridgeman, R. Canary and A. Sambarino, ``Bending, Entropy and proper affine actions of surface groups,'' in preparation.

\bibitem{BCY} M. Bridgeman, R. Canary and A.Yarmola, ``An improved bound for Sullivan's convex hull theorem,'' 
{\em Proc. LM.S.} {\bf 12}(2016), 146--168.

\bibitem{jp-neza} J.-P. Burelle and N. $\check{\mathrm{Z}}$ager Korenjak, ``Proper affine deformations of positive representations,'' preprint,
arXiv:2405.14658.

\bibitem{BIW-survey} M. Burger, A. Iozzi and A. Wienhard, ``Higher Teichm\"uller spaces: From $\mathsf{SL}(2,\mathbb R)$ to other Lie groups,''
{\em Handbook of Teichm\"uller Theory Vol.  IV,} (2014), European Mathematical Society, 539--618.

\bibitem{canary-tsouvalas} R. Canary and K. Tsouvalas, ``Topological restrictions on Anosov representations,''
{\em J. Topology} {\bf 13}(2020), 1497--1520. 

\bibitem{CZZ2} R. Canary, T. Zhang and A. Zimmer, ``Entropy rigidity for cusped Hitchin representations,'' preprint, arXiv:2201.04859.

\bibitem{CZZ3} R. Canary, T. Zhang and A. Zimmer, ``Patterson--Sullivan measures for transverse groups,'' {\em J. Mod. Dyn.} {\bf 20}(2024), 319--377.

\bibitem{CZZ4} R. Canary, T. Zhang and A. Zimmer, ``Patterson--Sullivan measures for relatively Anosov groups,'' 
{\em Math. Ann.}   {\bf 392}(2025), 2309--2363.

\bibitem{CZZ-survey} R. Canary, T. Zhang and A. Zimmer, ``Geometry and dynamics of transverse groups,''  preprint, arXiv:2502.07271.

\bibitem{CZZ5} R. Canary, T. Zhang and A. Zimmer, ``A rigidity theorem for complex Kleinian groups,'' preprint, arXiv:2511.20949.

\bibitem{choi-lee-marquis} S. Choi, G.-S. Lee and L. Marquis, ``Convex projective generalized Dehn filling,''
{\em Ann. E.N.S.} {\bf 53}(2020), 217--266.

\bibitem{chow} W.-L. Chow, ``On compact complex algebraic varieties,'' {\em Am. J. Math.} {\bf 71}(1949), 893--914.



\bibitem{DOP} F. Dal'bo, J.-P. Otal, and M. Peign\'e, ``S\'eries de Poincar\'e des groupes g\'eom\'etriquement finis,'' {\em Israel J. Math.} {\bf 118}(2000), 109--124.

\bibitem{DDGS} J. Danciger, T. Drumm, W. Goldman and I. Smilga, ``Proper actions of discrete groups of affine transformations,'' 
in {\em Dynamics, Geometry, Number Theory: the Impact of Margulis on Modern Mathematics}, University of Chicago Press (2022), 95--168.

\bibitem{DGK1} J. Danciger, F. Gu\'eritaud and F. Kassel, ``Geometry and topology of complete Lorentz spacetimes of constant curvature,''
{\em Ann. E.N.S.} {\bf 49}(2016), 1---56.

\bibitem{DGK2} J. Danciger, F. Gu\'eritaud and F. Kassel, ``Margulis spacetimes via the arc complex,''
{\em Invent. Math.} {\bf 204}(2016),  133--193.

\bibitem{DGK-aff} J. Danciger, F. Gu\'eritaud and F. Kassel, ``Proper affine actions of right-angled Coxeter groups,''
{\em Duke Math. J.} {\bf 169} 2020, 2231--2280.

\bibitem{DGK-CC} J. Danciger, F. Gu\'eritaud and F. Kassel, ``Convex cocompact actions in real projective
geometry,'' {\em Ann. E.N.S.}, to appear, arXiv:1704.08711.

\bibitem{DGK-combo}  J. Danciger, F. Gu\'eritaud and  F. Kassel, ``Combination theorems in convex projective geometry,'' preprint, arXiv:2407.09439.

\bibitem{DGKLM}  J. Danciger, F. Gu\'eritaud, F. Kassel, G,-S.Lee and L. Marquis,
``Convex cocompactness for Coxeter groups,'' {\em J.E.M.S.} {\bf 27}(2025), 119--181.

\bibitem{dey-borel} S. Dey, ``On Borel Anosov subgroups of $\mathsf{SL}(d, \mathbb R)$,'' {\em Geom. Top.} {\bf 29}(2025), 171--192.

\bibitem{DK-PS} S. Dey and M. Kapovich, ``Patterson-Sullivan theory for Anosov subgroups.,'' {\em Trans. A.M.S.} {\bf 375}(2022), 8687--8737.

\bibitem{DK-FP} S. Dey and M. Kapovich, ``Klein-Maskit combination theorem for Anosov subgroups: Free products,'' {\em Math. Z.} {\bf 305}(2023),
25 pp.

\bibitem{DK-KMA} S. Dey and M. Kapovich, ``Klein-Maskit combination theorem for Anosov subgroups: Amalgams,'' 
{\em Crelle Math.J} {\bf 819}(2025), 1--43.


\bibitem{DKL} S. Dey, M. Kapovich and B Leeb, ``A combination theorem for Anosov subgroups,''
{\em Math. Z.} {\bf 293}(2019),:551--578.
\bibitem{dey-tsouvalas} S. Dey and K. Tsouvalas, ``Anosov representations of amalgams,'' preprint, arXiv:2504.21802.

\bibitem{DKO} S. Dey, D. Kim,  and H. Oh, ``Ahlfors regularity of Patterson-Sullivan measures of Anosov groups and applications,''
preprint, arXiv:2401.12398.

\bibitem{DFWZ} S. Douba, B. Fl\'echelles, T. Weisman and F, Zhu, ``Cubulated hyperbolic groups admit Anosov representations,''
{\em Geom. Top.}, to appear, arXiv:2309.03695.

\bibitem{douba-tsouvalas} S. Douba and K. Tsouvalas, ``Anosov groups that are indiscrete in rank one,'' 
{\em I.M.R.N.} {\bf 2024}, 8377--8383.

\bibitem{dular-schlenker} B. Dular and J.M. Schlenker, ``Convex co-compact hyperbolic manifolds are determined by their pleating lamination,''
preprint, arXiv:2403.10090.

\bibitem{drumm} T. Drumm,  ``Fundamental polyhedra for Margulis space-times", {\em  Topology} {\bf 31}(1992), 677--683.

\bibitem{DG-crooked} Drumm, T. and Goldman, W., ``The geometry of crooked planes," {\em Topology},  {\bf 38}(1999), 323--351.

\bibitem{epstein-marden} D.B.A. Epstein and A. Marden, ``Convex hulls in hyperbolic space, a theorem of Sullivan, and measured pleated surfaces" {\em Analytical and Geometrical Aspects of Hyperbolic Space}, Cambridge University Press, 1987, 113--253.

\bibitem{FPV} J.Farre, M.B. Pozzetti, and G. Viaggi, ``Topological and geometric restrictions on hyperconvex representations,'' 
{\em J. Reine Angew. Math}, to appear,  arXiv:2403.13668.

\bibitem{FPV2} J.Farre, M.B. Pozzetti, and G. Viaggi, ``Geometry of hyperconvex representations of surface groups,'' preprint, arXiv:2407.20071.

\bibitem{fried-goldman} D. Fried and W. Goldman, ``Three-dimensional affine crystallographic groups,''
{\em Adv. in Math.} {\bf 47}(1983),  1--49.

\bibitem{ghosh}
S. Ghosh, ``Deformation of fuchsian representations and proper affine actions,'' preprint,
arXiv:2312.16655.

\bibitem{ghosh-treib} S. Ghosh and N. Treib, ``Affine Anosov representations and proper actions,'' {\em IMRN}  {\bf 16}(2023), 14334--14367.

\bibitem{GMT} O. Glorieux, D. Montclair and N. Tholozan, ``Hausdorff dimension of limit sets for projective Anosov representations,''
{\em Jour. l'\'Ecole Poly.} {\bf 10}(2023), 1157--1193.


\bibitem{GLM} W.  Goldman, F. Labourie, G. A. Margulis, ``Proper affine actions and geodesic flows of hyperbolic surfaces", 
{\em Ann. of Math.} {\bf 170}, 2009,  1051--1083.

\bibitem{GGKW} F. Gu\'eritaud, O. Guichard, F. Kassel and A. Wienhard, ``Anosov representations and proper actions,''
{\em Geom. Top.} {\bf 21}(2017), 485--584.

\bibitem{GW} O. Guichard and A. Wienhard, ``Anosov representations: Domains of discontinuity and
applications,'' {\em Invent. Math.} {\bf 190}(2012), 357--438.

\bibitem{haissinsky} P. Haissinsky, ``Hyperbolic groups with planar boundaries,''  {\em Invent. Math.} {\bf  201}(2015),  239--307.

\bibitem{hironaka} H. Hironaka, ``Resolution of singularities of an algebraic variety over a field of characteristic zero. {I}, {II}," 
{\em Ann. of Math.} {\bf 79}(1964), 109--203.

\bibitem{IZ} M. Islam and A. Zimmer, ``Convex co-compact actions of relatively hyperbolic groups,''
{\em Geom. Top.} {\bf 27}(2023), 417--511.

\bibitem{IZ2} M. Islam and A. Zimmer, ``Convex co-compact groups with one dimensional boundary faces,''
{\em Grps. Gyn. Dyn.} {\bf 18}(2024), 1145--1183.

\bibitem{KLM} J. Kahn, F. Labourie and S. Mozes, ``Surface subgroups in uniform lattices of some semi-simple Lie groups,''
{\em Acta Math.} {\bf 232}(2024),  79--220.

\bibitem{kahn-markovic} J. Kahn and V. Markovic, ``Immersing almost geodesic surfaces in a closed hyperbolic 3-manifold,''
{\em Ann. of Math.} {\bf 175}(2012),  1127--1190 

\bibitem{kapovich-benakli} I. Kapovich and N. Benakli, ``Boundaries of hyperbolic groups,'' in {\em Combinatorial and Geometric
Group Theory}, Contemp. Math. vol. 296(2002), A.M.S., 39--93.

\bibitem{kapovich-gt} M. Kapovich, ``Convex projective structures on Gromov-Thurston manifolds,'' {\em Geom. Top.} {\bf 11}(2007)  1777--1830. 

\bibitem{KL-survey} M. Kapovich and B. Leeb, ``Discrete isometry groups of symmetric spaces,'' {\em Handbook of Group Actions vol. IV}. 
International Press (2018), 191--290. 

\bibitem{KL} M. Kapovich and B. Leeb, ``Relativizing characterizations of Anosov subgroups I,'' {\em  Groups, Geom. Dyn.} {\bf 17}(2023)  1005--1071.

\bibitem{KLP-survey} M. Kapovich, B. Leeb and J. Porti, ``Some recent results on Anosov representations,'' 
{\em Trans. Groups} {\bf 21}(2016), 1105--1121. 

\bibitem{KLP} M. Kapovich, B. Leeb and J. Porti, ``Anosov subgroups: Dynamical and geometric characterizations,''
{\em Eur. Math. J.} {\bf 3}(2017), 808--898.



\bibitem{kassel-icm} F. Kassel, ``Geometric structures and representations of discrete groups,''
in {\em Proceedings of the ICM 2018, vol. 2,} World Scientific (2019), 1113--1150.

\bibitem{kassel-survey} F. Kassel, ``Discrete subgroups of semisimple Lie groups, beyond lattices,'' in
{\em Groups St Andrews 2022 in Newcastle}, LMS Lecture Notes  vol. 496 (2024), Cambridge University Press, 118--190.

\bibitem{kassel-smilga}
F. Kassel and I. Smilga, ``Affine properness criterion for $\mathrm{A}$nosov representations and generalizations,'' in preparation.


\bibitem{KO-RA} D. Kim and H. Oh, ``Relatively Anosov groups: finiteness, measure of maximal entropy and reparameterization,''
{\em Crelle's J.}, to appear, arXiv:2404.09475. 

\bibitem{KOW} D. Kim, H. Oh and Y. Wang, ``Properly discontinuous actions, growth indicators and conformal measures for transverse subgroups," preprint, arXiv:2306.06846.

\bibitem{klein} F.Klein,  ``Neue Beitr\"age zur Riemann'schen Functionentheorie,'' {\em Math. Ann} {\bf 21}(1883), 141--218.

\bibitem{kouron} C. Kourouniotis, ``The geometry of bending quasi-Fuchsian groups,'' in {\em Discrete Groups and Geometry (Birmingham, 1991)}
 L.M.S. Lecture Note Series vol. 74, 1992, 148--164.


\bibitem{kouron-cont} C. Kourouniotis, ``On the continuity of bending''
{\em Geometry and Topology Monographs}, {\bf Vol 1}, (1998) 317--334.

\bibitem{labourie-invent} F. Labourie, ``Anosov flows, surface groups and curves in projective space,'' 
{\em Invent. Math.} {\bf 165} (2006), no.~1, 51--114.

\bibitem{lahn1} M. Lahn, ``Reducible suspensions of Anosov representations,'' {\em Grps. Geom. Dyn.}, to appear, arXiv:2312.09886.

\bibitem{lahn2} M. Lahn, ``Which reducible representations are Anosov?,'' preprint, arXiv:2411.15321.

\bibitem{MMMZ1} S. Maloni, G. Martone,  F. Mazzoli and T. Zhang,  "d-pleated surfaces and their shear-bend coordinates,"
preprint,  arXiv:2305.11780.

\bibitem{MMMZ2} S. Maloni, G. Martone,  F. Mazzoli and T. Zhang, ``Topology of the space of $d$-pleated surface,''
preprint,  arXiv:2508.04813.

\bibitem{Mar1} G. Margulis, ``Free properly discontinuous groups of affine transformations", 
{\em Dokl. Akad.Nauk SSSR} {\bf 272}, 1983, 785--788. 

\bibitem{Mar2} G. Margulis, ``Complete affine locally flat manifolds with a free fundamental group", {\em J.
Soviet Math.} {\bf 134}, 1987, 129--134.

\bibitem{markovic} V.  Markovic,  ``Quasisymmetric groups,''  {\em J. Amer. Math. Soc.} {\bf 19}(2006),  673--715.

\bibitem{marquis-survey} L. Marquis, ``Around groups in Hilbert geometry,'' in
{\em Handbook of Hilbert geometry}, European Mathematical Society (2014), 207--261. 

\bibitem{maskit} B. Maskit, ``On Klein's combination theorem III,'' in {\em Advances in the theory of Riemann surfaces},
Princeton University Press (1971), 297--316.

\bibitem{MV} F. Mazzoli and G. Viaggi, ``$\mathsf{SO}_0(2,n+1)$-maximal representations and hyperbolic surfaces,''
{\em Mem. A.M.S.} {\bf 313}(2025).

\bibitem{MP} P. Menal-Ferrer and J. Porti, ``Twisted cohomology for hyperbolic 3-manifolds,'' 
{\em Osaka J. Math.} {\bf 49}(2012),  741--769.

\bibitem{Mess} G. Mess, ``Lorentz spacetimes of constant curvature", {\em Geometriae Dedicata}, {\bf 126}, 
2007, 3--45

\bibitem{milnor} J.Milnor, ``On fundamental groups of complete affinely flat manifolds,''
{\em Adv. Math.} {\bf 25}(1977), 178--187.

\bibitem{morgan} J. Morgan, ``On Thurston's uniformization theorem for three-dimensional manifolds'' in 
{\em The Smith Conjecture}, Pure Appl. Math. {\bf 112}, Academic Press, 1984, 37--125.

\bibitem{namazi-souto} H. Namazi and J. Souto, ``Non-realizability and ending
laminations,'' {\em Acta Math.} {\bf 209}(2012),  323--395.

\bibitem{ohshika-density} K. Ohshika, ``Realising end invariants by limits of minimally parabolic, geometrically finite
groups,'' {\em Geom. Top.} {\bf 15}(2011), 827--890.

\bibitem{patterson} S.J. Patterson, ``The limit set of a Fuchsian group,'' {\em Acta Math.} {\bf 136}(1976), 241--273. 

\bibitem{porti} J. Porti, ``Dimension of representation and character varieties for two- and three-orbifolds,''
{\em Alg. Geom. Top.} {\bf 22}(2022), 1905--1967.



\bibitem{PSW1} M. Pozzetti, A. Sambarino and A. Wienhard, ``Conformality for a robust class of non-conformal attractors,''
{\em J. Reine Angew. Math} {\bf 74}(2021), 1--51.

\bibitem{quint} J.F. Quint, ``Mesures de Patterson-Sullivan en rang sup\'erieur,'' {\em G.A.F.A.} {\bf 12}(2002), 776--809.

\bibitem{ruelle} D. Ruelle, ``Repellers for real analytic maps,''
{\em Ergodic Theory Dynamical Systems} {\bf 2}(1982),  99--107. 

\bibitem{sambarino-quant} A. Sambarino, ``Quantitative properties of convex representations,''
{\em Comm. Math. Helv.} {\bf 89}(2014), 443--488.

\bibitem{sambarino-dichotomy} A. Sambarino, ``A report on an ergodic dichotomy,'' {\em Ergod. Theory Dyn. Syst.} (2024), 236--289. 

\bibitem{sambarino} A. Sambarino, ``Asymptotic properties of infinitesmal characters and applications,'' preprint, arXiv:2406.06250.

\bibitem{schwartz-pappus} R. Schwartz, ``Pappus's Theorem and the Modular Group,''
{\em I.H.E.S. Pub. Math.} {\bf 78}(1993), 187--206.

\bibitem{schwartz-CR} R. Schwartz, {\em Spherical CR Geometry and Dehn Surgery}
Princeton University Press, 2007.

\bibitem{schwartz-anosov} R. Schwartz, ``Patterns of Geodesics, Shearing, and Anosov Representations of the Modular Group,''
preprint, arXiv:2412.18457.

\bibitem{smilga-ortho} I. Smilga, ``Fundamental domains for properly discontinuous affine groups,''
{\em Geom. Ded.}  {\bf 171}(2014), 203--229.

\bibitem{Smig} I. Smilga, ``Proper affine actions on semisimple Lie algebras", {\em Ann. Inst. Fourier
(Grenoble)}, {\bf 66}(2016), 785--831

\bibitem{sullivan-ergodic} D. Sullivan, ``The ergodic theory at infinity of an arbitrary discrete group of hyperbolic motions,'' in
{\em Riemann surfaces and related topics: proceedings of the 1978 Stony Brook conference}, Princeton University
Press (1981), 465--496.

\bibitem{sullivan-hd} D. Sullivan, ``Entropy, Hausdorff measures old and new, and limit sets of geometrically finite Kleinian groups,''
{\em Acta Math.} {\bf153(}1984), 259--277.

\bibitem{tholozan-tsouvalas} N. Tholozan and K. Tsouvalas, ``Linearity and indiscreteness of amalgamated products of hyperbolic groups,''
{\em I.M.N.R.} {\bf 2023}, 21290--21319.

\bibitem{tomanov} G. Tomanov, ``Properly discontinuous group actions on affine homogeneous spaces,''
{\em Proc. Steklov Inst. Math.} {\bf 292}(2016),  260--271.

\bibitem{thurston-notes} W. Thurston, {\em The geometry and topology of 3-manifolds},
lecture notes, Princeton University, 1980, available at: \url{http://library.msri.org/nonmsri/gt3m/}

\bibitem{traaseth-weisman} A. Traaseth and T. Weisman, ``Combination theorems for geometrically finite convergence groups,''
{\em Alg. Geom. Top.}, to appear, arXiv:2305.08011.

\bibitem{tsouvalas-borel} K. Tsouvalas, ``On Borel Anosov representations in even dimensions,''
{\em Comm. Math. Helv.} {\bf 95}(2020) 749--763.

\bibitem{tsouvalas-density} K. Tsouvalas, ``Robust quasi-isometric embeddings inapproximable by Anosov representations,''
preprint, arXiv:2402.09339.   

\bibitem{tsouvalas-sep} K. Tsouvalas, ``Separability of unipotent-free abelian subgroups in linear groups,'' preprint, arXiv:2504.19995.

\bibitem{tsouvalas-weisman} K. Tsouvalas and T. Weisman, ``Singular value gap estimates for free products of semigroups,'' preprint,
arXiv:2409.20330.

\bibitem{weisman-EGF} T. Weisman, ``An extended definition of Anosov representation for relatively hyperbolic groups,'' preprint,  arXiv:2205.07183.

\bibitem{weisman-EGF2} T. Weisman, ``Examples of extended geometrically finite representations,'' preprint,  arXiv:2311.18653.

\bibitem{weisman-dehn} T. Weisman, ``Dehn filling in semisimple Lie groups,'' preprint, arXiv: 2502.17592.

\bibitem{wienhard-icm}  A. Wienhard, ``An invitation to higher Teichm\"uller theory,'' in {\em  Proceedings of the ICM 2018 vol. 2}, 
World Scientific (2019), 1013--1039.

\bibitem{neza} N. $\check{\mathrm{Z}}$ager Korenjak, ``Constructing proper affine actions via higher strip deformations,'' preprint,
arXiv:2205.14712.

\bibitem{zhang-zimmer}  T. Zhang and A. Zimmer, ``Regularity of limit sets of Anosov representations,'' {\em J. Topology} {\bf 17}(2024), 72 pp.

\bibitem{zhu-thesis} F. Zhu, ``Relatively dominated representations,''
{\em Ann. Inst. Fourier} {\bf 71}(2021),  2169--2235.

\bibitem{ZZ1} F. Zhu and A. Zimmer, ``Relatively Anosov representations via flows I: theory,'' {\em Grps. Geom, Dyn.}, to appear, arXiv:2207.14737.

\bibitem{ZZ2} F. Zhu and A. Zimmer, ``Relatively Anosov representations via flows II: examples,'' {\em  Jour. L.M.S.}  {\bf 109}(2024), 61 pages.

\bibitem{zimmer-cc} A. Zimmer, ``Projective Anosov representations, convex cocompact actions, and rigidity,''
{\em J. Diff. Geom.} {\bf 119}(2021),  513--586.

\end{scriptsize}
\end{thebibliography}
\end{document}